%% file: thesismaster.tex
\documentclass[defaultstyle,11pt]{thesis}

\usepackage{graphics, color}
\usepackage{amsfonts}
\usepackage{amssymb}
\usepackage{latexsym}

\title{Two Gauss-Bonnet and Poincar\'{e}-Hopf Theorems for Orbifolds with Boundary}

\author{Christopher W.}{Seaton}

\otherdegrees{B.A., Kalamazoo College, 1999}

\degree{Doctor of Philosophy}
	{Ph.D., Mathematics}
	
\dept{Department of}
	{Mathematics}

\advisor{Assoc. Prof.}
	{Carla Farsi}

\reader{Siye Wu}									%  2nd person to sign thesis
\abstract{  \OnePageChapter
	The goal of this work is to generalize the Gauss-Bonnet and Poincar\'{e}-Hopf Theorems
	to the case of orbifolds with boundary.  We present two such generalizations, the first
	in the spirit of \cite{satake2}, which uses an argument parallel to that contained in
	\cite{sha}.  In this case, the local data (i.e. integral of the curvature in the case
	of the Gauss-Bonnet Theorem and the index of the vector field in the case of the
	Poincar\'{e}-Hopf Theorem) is related to Satake's orbifold Euler characteristic, a
	rational number which depends on the orbifold structure.

	For the second pair of generalizations, we use a more recent orbifold cohomology
	\cite{chenruan} to express the local data in a way which can be related to the Euler
	characteristic of the underlying space of the orbifold.  This case applies only to
	orbifolds which admit almost-complex structures.
	}

	\dedication[Dedication]{					% NEVER use \OnePageChapter here.
	To Michael Thomas Seaton and Christina Heather Bost, my siblings old and new,
	with congratulations and best wishes.
		}

	\acknowledgements{	\OnePageChapter		% *MUST* BE ONLY ONE PAGE!

	I am honored to take this opportunity to thank my advisor, Carla Farsi, for her encouragement,
	instruction, and inspiration throughout my time at Colorado.  This work would not have been
	possible without her.

	I would also like to thank Alexander Gorokhovsky, Erich McAlister, Arlan Ramsay,
	Yongbin Ruan, and Siye Wu for many helpful discussions and suggestions,
	and John Massman for helpful discussions and pointers about using \LaTeX .  
	As well, thanks to Lynne Walling for encouragement.

	Thanks are also due to the kids who would answer the phone at crazy hours
	of the night to listen to my frustration.  To
	L.J., Demangy, Andrea, Crickey, Her Majesty, and my brother, thank you for your
	faith, support, and patience.

	And of course I would like to thank my parents for making everything possible.
	}

\ToCisShort							%  a 1-page Table of Contents

\LoFisShort							%  a 1-page List of Figures
\emptyLoT							%  no List of Tables at all
\begin{document}

\input thesismacros.tex					%  macro definitions
\input thesisch1.tex					% Chapter 1 body
\input thesisch2.tex					% Chapter 2 body
\input thesisch3.tex					% Chapter 3 body
\input thesisch4.tex					% Chapter 4 body

% XXXXXXXXXXXXXXXXXXXXXXXXXXXXXXXXXXXXXXXXXXXXXXXXXXXXXXXXXXXXXX
% XXXXXXXXXX			bibliography		XXXXXXXX
% XXXXXXXXXX	commented out until i fix refs format	XXXXXXXX
% XXXXXXXXXXXXXXXXXXXXXXXXXXXXXXXXXXXXXXXXXXXXXXXXXXXXXXXXXXXXXX

\bibliographystyle{plain}					%  or "siam", or "alpha", or "abbrv"
								%  see other styles (.bst files) in
								%  $TEXHOME/texmf/bibtex/bst

\bibliography{thesisrefs}					%  bib database file thesisrefs.bib

% XXXXXXXXXXXXXXXXXXXXXXXXXXXXXXXXXXXXXXXXXXXXXXXXXXXXXXXXXXXXXX
% XXXXXXXXXX			appendices		XXXXXXXX
% XXXXXXXXXXXXXXXXXXXXXXXXXXXXXXXXXXXXXXXXXXXXXXXXXXXXXXXXXXXXXX

 \appendix	% don't forget this line if you have appendices!

\input thesisappA.tex			% file with Appendix A contents
 \input thesisappB.tex			% file with Appendix B contents

% XXXXXXXXXXXXXXXXXXXXXXXXXXXXXXXXXXXXXXXXXXXXXXXXXXXXXXXXXXXXXX
% XXXXXXXXXX			the end			XXXXXXXX
% XXXXXXXXXXXXXXXXXXXXXXXXXXXXXXXXXXXXXXXXXXXXXXXXXXXXXXXXXXXXXX

\end{document}

%% file: thesismacros.tex
%	Macros			Christopher Seaton. seatonc@colorado.edu

	%%%%%	DEFINE RATIONAL, REAL, AND COMPLEX NUMBER SIGNS
\newcommand{\C}{\mathbb{C}}
\newcommand{\Q}{\mathbb{Q}}
\newcommand{\R}{\mathbb{R}}
\newcommand{\Z}{\mathbb{Z}}

	%%%%%	SET UP THEOREM ENVIRONMENT
\newtheorem{theorem}{Theorem}[section]

	%%%%%	SET UP DEFINITION ENVIRONMENT	
\newtheorem{definition}[theorem]{Definition}

	%%%%%	SET UP LEMMA ENVIRONMENT	
\newtheorem{lemma}[theorem]{Lemma}

	%%%%%	SET UP EXAMPLE ENVIRONMENT	
\newtheorem{example}[theorem]{Example}

	%%%%%	SET UP PROPOSITION ENVIRONMENT
\newtheorem{proposition}[theorem]{Proposition}

	%%%%%	SET UP CLAIM ENVIRONMENT
\newtheorem{claim}[theorem]{Claim}

	%%%%%	SET UP COROLLARY ENVIRONMENT
\newtheorem{corollary}[theorem]{Corollary}

	%%%%%	SET UP REMARK ENVIRONMENT
\newtheorem{remark}[theorem]{Remark}

%% file: thesisch1.tex
% Thesis Chapter 1, Christopher Seaton, seatonc@colorado.edu
% XXXXXXXXXXXXXXXXXXXXXXXXXXXXXXXXXXXXXXXXXXXXXXXXXXXXXXXXXXXXXXXXXXX
% XXXXXXXXXXXXXXXXXXXXXXXXXXXXXXXXXXXXXXXXXXXXXXXXXXXXXXXXXXXXXXXXXXX
%			Chapter 1: Introduction
% XXXXXXXXXXXXXXXXXXXXXXXXXXXXXXXXXXXXXXXXXXXXXXXXXXXXXXXXXXXXXXXXXXX
% XXXXXXXXXXXXXXXXXXXXXXXXXXXXXXXXXXXXXXXXXXXXXXXXXXXXXXXXXXXXXXXXXXX
%
%	should just be one short section
%

\chapter{Introduction}

An orbifold is perhaps the simplest case of a singular space;
it is a topological space which is locally diffeomorphic to
$\R^n/G$ where $G$ is a finite group.  Orbifolds were originally
introduced by Satake in \cite{satake1} and \cite{satake2}, where
they were given the name $V$-manifold, and rediscovered by Thurston
in \cite{thurston}, where the term orbifold was coined.  Satake
and Thurston's definitions differ, however, in that Satake required
the group action to have a fixed point set of codimension at least
two, while Thurston did not.  Hence, Thurston's definition allows group
actions such as reflections through hyperplanes.
Today, authors differ on whether or not this requirement is made;
often, when it is, the orbifolds are referred to as {\bf codimension-2
orbifolds}.  It is these orbifolds which are our object of study.

The point of view of this work is that an orbifold structure is a
generalization of a differentiable structure on a manifold.  We do not
mean to suggest that the underlying space of an orbifold is necessarily
a topological manifold; this is only the case in dimension $\leq 2$, and
not even in dimension 1 if the codimension-2 requirement is lifted.  However,
there are many examples of orbifolds whose underlying topological spaces
are indeed manifolds.  In these cases, we view the orbifold structure as
a singular differentiable structure on the manifold.  It should be noted
that this is used as a guiding principle only, and that our results
apply in general.

Hence, we improve upon Satake's Gauss-Bonnet theorem for
orbifolds \cite{satake2} by developing a Gauss-Bonnet integrand (and
corresponding {\bf orbifold Euler Class}) whose integral relates to
the Euler Characteristic of the underlying topological space, as opposed
to the orbifold
Euler Characteristic (Theorem \ref{gbcr}).  This result depends on recent developments in the
theory of orbifolds, most notably the Orbifold Cohomology Theory of
Chen-Ruan \cite{chenruan}, and hence is restricted to the case of an
orbifold which admits an almost complex structure.

As is well-known, the Gauss-Bonnet theorem is very closely related to the
Poincar\'{e}-Hopf theorem; indeed, either of the two theorems can be
viewed as a corollary of the other.  With a new Gauss-Bonnet theorem,
then, comes a new Poincar\'{e}-Hopf theorem.  Following the work of
Sha on the secondary Chern-Euler class for manifolds \cite{sha}, we
generalize the Poincar\'{e}-Hopf theorem to the case of orbifolds with
boundary (Corollary \ref{phcr}).

To this end, a remark is in order.  Satake's original definition of an orbifold with boundary
was in many senses not very strict.  In particular, the boundaries of his
orbifolds were not necessarily orbifolds.
This allows little control over vector fields on the boundary; indeed,
it is not generally the case that the orbifold is locally a product near the
boundary, and hence there need not exist a vector field which does not vanish
on the boundary.  Since then, a more natural
definition of orbifold with boundary has been given in \cite{thurston}.
Using this definition, we are able to refine Satake's original Gauss-Bonnet
theorem for the case with boundary using his original techniques.

The outline of this work is roughly as follows.  In Chapter \ref{chbackground},
we collect the necessary background information on orbifolds and orbifolds
with boundary, including several examples, paying particular attention
to the behavior of vector fields on orbifolds.  In Chapter \ref{chgbph}, we
review Satake's Gauss-Bonnet and Poincar\'{e}-Hopf theorems for orbifolds
and orbifolds with boundary, making improvements where possible using the
modern definition of an orbifold with boundary (see Theorem \ref{gbb1}).
We also apply the
arguments of Sha \cite{sha} to characterize the boundary term in the
case with boundary as the evaluation of a secondary characteristic class on the boundary (see Theorem \ref{pcb1}).
It is easy to see that, even in the
case of a manifold, the boundary term of this formula will always depend
on the vector field.  For instance, consider a fixed vector field on $S^2$
with one singular point $p$.  Remove an open disk to produce
a vector field on a manifold with boundary.  The index of the vector field
depends on whether $p$ is contained in the disk removed, but the Euler
Characteristic does not; hence, the boundary term must depend on the
vector field.  The spirit of the Poincar\'{e}-Hopf theorem is that this
term should be formulated in a manner as independent of the vector field as
possible; it is this reason that we chose the result of Sha to generalize
to orbifolds.

In Chapter \ref{chcr}, we review the Chen-Ruan orbifold cohomology, and extend
it in a straightforward manner to the case with boundary.  Loosely speaking,
the idea of this cohomology theory is to associate to an orbifold $Q$
another orbifold, $\tilde{Q}$ (where at least one of the connected
components of $\tilde{Q}$ is diffeomorphic to $Q$), and use the cohomology
groups of $\tilde{Q}$ (we should note that this informal description leaves
out an important modification to the grading of the groups).  Using this cohomology
theory, we develop an Euler Class which relates to the Euler Characteristic
of the underlying topological space of $Q$.  The essential idea here is to
apply the Chern-Weil description of characteristic classes to the curvature
of a connection on $\tilde{Q}$, yielding a (non-homogeneous) characteristic
class in orbifold cohomology.  Similarly, the index of a vector field $X$
on $Q$ is computed to be the index of its pull-back $\tilde{X}$ onto
$\tilde{Q}$.  This suggests the paradigm that geometric structures on
$\tilde{Q}$ can be considered to be structures on $Q$ which take multiple
values on singular sets.  It is in this manner that we prove the
aforementioned Gauss-Bonnet and Poincar\'{e}-Hopf theorems for orbifolds
(Theorem \ref{gbcr} and Corollary \ref{phcr})
and orbifolds with boundary
(Theorems \ref{gbbdcr} and \ref{pccr}), relating to the Euler Characteristics of
the underlying space.

Theorems, definitions, examples, etc. are numbered sequentially according to the
section in which they appear.  So Theorem $X.Y.Z$ is in Chapter $X$, Section $Y$,
and follows Definition (or Lemma, Example, etc.) $X.Y.Z-1$.  Figures are
numbered independently according to the Chapter in which they appear.

%% file: thesisch2.tex
% Thesis Chapter 2, Christopher Seaton, seatonc@colorado.edu
% XXXXXXXXXXXXXXXXXXXXXXXXXXXXXXXXXXXXXXXXXXXXXXXXXXXXXXXXXXXXXXXXXXX
% XXXXXXXXXXXXXXXXXXXXXXXXXXXXXXXXXXXXXXXXXXXXXXXXXXXXXXXXXXXXXXXXXXX
%			Chapter 2: Orbifolds and Their Structure
% XXXXXXXXXXXXXXXXXXXXXXXXXXXXXXXXXXXXXXXXXXXXXXXXXXXXXXXXXXXXXXXXXXX
% XXXXXXXXXXXXXXXXXXXXXXXXXXXXXXXXXXXXXXXXXXXXXXXXXXXXXXXXXXXXXXXXXXX
%
%	Subsection 1: Definitions and Examples
%					orbifolds and orbifolds with boundary
%					examples
%					structures on orbifolds
% Subsection 2:	The Dimension of a Singularity
%					motivation
%					properties
%					examples of 0-dim singularities
%

\chapter{Orbifolds and Their Structure}
\label{chbackground}

% XXXXXXXXXXXXXXXXXXXXXXXXXXXXXXXXXXXXXXXXXXXXXXXXXXXXXXXXXXXXXXXXXXX
%			Section: Definitions and Examples
% XXXXXXXXXXXXXXXXXXXXXXXXXXXXXXXXXXXXXXXXXXXXXXXXXXXXXXXXXXXXXXXXXXX

\section{Definitions and Examples}

In this section, we collect the definitions and background we will need.
For more information, the reader is referred to the original work of Satake in
\cite{satake1} and \cite{satake2}.  As well,
\cite{ruangwt} contains as an appendix a thorough introduction to orbifolds,
focusing on their differential geometry.  Other good introductions include
\cite{thurston} and \cite{borzellino}, the former providing a great deal of
information on the topology of low-dimensional orbifolds.  The reader is warned that the
definition used in these latter two works is more general than ours, as it admits
group actions which fix sets of codimension 1.  For the most part, we follow
the spirit of Satake and Ruan.

% XXXXXXXXXXXXXXXXXXXXXXXXXXXX	SUBSUBSECTION: ORBIFOLDS AND ORBIFOLDS WITH BOUNDARY

\subsection{Orbifolds and Orbifolds With Boundary}
\label{orbandwbd}

Let $X_Q$ be a Hausdorff space.

% XXXXXXXXXXXXXXXXXXXXXXXXXXXX	DEFINITION

\begin{definition}[orbifold chart]
Let $U \subset X_Q$ be a connected open set.
A (${\mathcal C}^\infty$) {\bf orbifold chart} for $U$ (also known as a
(${\mathcal C}^\infty$) {\bf local unifomizing system})
is a triple $\{ V, G, \pi \}$ where
\begin{itemize}
\item		$V$ is an open subset of $\R^n$,
\item		$G$ is a finite group with a (${\mathcal C}^\infty$) action on $V$ such
		that the fixed point
		set of any $\gamma \in G$ which does not act trivially on $V$ has
		codimension at least 2 in $V$, and
\item		$\pi : V \rightarrow U$ is a surjective continuous map such that
		$\forall \gamma \in G$, $\pi \circ \gamma = \pi$ that induces a
		homeomorphism $\tilde{\pi} : V/G \rightarrow U$.
\end{itemize}
If $G$ acts effectively on $V$, then the chart is said to be
{\bf reduced}.
\end{definition}

The definition of the appropriate notion of `chart' for orbifolds with boundary is
similar:

\begin{definition}[orbifold chart with boundary]
Let $U \subset X_Q$ be a connected open set.
A (${\mathcal C}^\infty$) {\bf orbifold chart with boundary}
or (${\mathcal C}^\infty$) {\bf local unifomizing system}
for $U$ is a triple $\{ V, G, \pi \}$ where
\begin{itemize}
\item		$V$ is an open subset of
				$\R_+^n : = \{ (x_1, x_2, \ldots \, x_n) \in \R^n : x_1 \geq 0 \}$,
\item		$G$ is a finite group with a (${\mathcal C}^\infty$) action on $V$ such
				that the fixed point
				set of any $\gamma \in G$ which does not act trivially on $V$ has
				codimension at least 2 in $V$, and such that
				$\gamma \partial \R_+^n \subset \partial \R_+^n$, and
\item		$\pi : V \rightarrow U$ is a surjective continuous map such that
				$\forall \gamma \in G$, $\pi \circ \gamma = \pi$ that induces a
				homeomorphism $\tilde{\pi} : V/G \rightarrow U$.
\end{itemize}
Again, if $G$ acts effectively on $V$, then the chart is said to be {\bf reduced}.
\end{definition}

If $V \cap \partial \R_+^n = \emptyset$, then $\{ V, G, \pi \}$ is an ordinary orbifold
chart; for emphasis, we may refer to these as {\bf orbifold charts without boundary}.
Note also that if $\{ V, G, \pi \}$ is an orbifold chart with boundary for some set
$U$, then restricting the chart to $\partial V = \partial \R_+^n \cap V$, it is clear that
$\{ \partial V , G , \pi_{| \partial V} \}$ is an orbifold chart without boundary
for $\pi (\partial V)$.

We will always use the notation that if $V_i$ is the domain of a chart, then $G_i$
is the group of the chart, $\pi_i$ the projection for the chart, and $U_i$ the
range of the chart; i.e. items with the same subscript correspond to the same chart.

Orbifold charts relate to one another via {\bf injections}.

% XXXXXXXXXXXXXXXXXXXXXXXXXXXX	DEFINITION

\begin{definition}[injection]
If $\{ V_i , G_i , \pi_i \}$ and $\{ V_j , G_j , \pi_j \}$ are two orbifold
charts (with or without boundary) for $U_i$ and $U_j$, respectively, where
$U_i \subset U_j \subset X_Q$, then an {\bf injection}
$\lambda_{ij} : \{ V_i , G_i , \pi_i \} \rightarrow \{ V_j , G_j , \pi_j \}$
is a pair $\{ f_{ij}, \phi_{ij} \}$ where
\begin{itemize}
\item		$f_{ij} : G_i \rightarrow G_j$ is an injective homomorphism such
		that if $K_i$ and $K_j$ denote the kernel of the action of $G_i$ and $G_j$,
		respectively, then $f_{ij}$ restricts to an isomorphism of $K_i$ onto $K_j$,
		and
\item		$\phi_{ij} : V_i \rightarrow V_j$ is a smooth embedding such that
		$\pi_i = \pi_j \circ \phi_{ij}$ and such that for each $\gamma \in G_i$,
		$\phi_{ij} \circ \gamma = f_{ij}(\gamma) \circ \phi_{ij}$.
		If the charts have boundary then $\phi_{ij}(\partial V_i) \subseteq \partial V_j$.
\end{itemize}
\end{definition}

Given an orbifold chart $\{ V, G, \pi \}$ (with or without boundary),
each $\gamma \in G$ induces an
injection $\lambda_\gamma$ of $\{ V, G, \pi \}$ into itself via
\[
\begin{array}{rcl}
	f_\gamma	&:&	G			\rightarrow		G			\\\\
			&:&	\gamma^\prime		\mapsto	
					\gamma \gamma^\prime \gamma^{-1}			\\\\
	\phi_\gamma	&:&	V			\rightarrow		V			\\\\
			&:&	x			\mapsto	\gamma x.
\end{array}
\]
Note that this injection is trivial if $\gamma$ acts trivially.  Similarly, given an injection
$\lambda_{ij} : \{ V_i , G_i , \pi_i \} \rightarrow \{ V_j , G_j , \pi_j \}$,
each element $\gamma \in G_j$ defines an injection
$\gamma \lambda_{ij} : \{ V_i , G_i , \pi_i \} \rightarrow \{ V_j , G_j , \pi_j \}$
with $\gamma \lambda_{ij} := \{ \gamma f_{ij}\gamma^{-1}, \gamma \circ \phi_{ij} \}$.
Moreover, every two injections of $\{ V_i , G_i , \pi_i \}$ into $\{ V_j , G_j , \pi_j \}$ are related in this manner (see \cite{satake2}, Lemma 1).

Two orbifold charts $\{ V_i, G_i, \pi_i \}$ and $\{ V_j, G_j, \pi_j \}$ are said to be
{\bf equivalent} if $U_i = U_j$, and there is an injection $\lambda_{ij}$ with
$f_{ij}$ an isomorphism and $\phi_{ij}$ a diffeomorphism.

% XXXXXXXXXXXXXXXXXXXXXXXXXXXX	DEFINITION

\begin{definition}[orbifold]
An {\bf orbifold} $Q$ is a Hausdorff space $X_Q$, the {\bf underlying space of $Q$},
together with a family ${\mathcal F}$ of orbifold charts (without boundary) such that
\begin{itemize}
\item		Each $p \in X_Q$ is contained in an open set $U_i$ covered by an orbifold
		chart $\{ V_i, G_i, \pi_i \} \in {\mathcal F}$.  If $p \in U_i \cap U_j$
		for $U_i$ and $U_j$ uniformized sets, then there is a uniformized set
		$U_k$ such that $p \in U_k \subset U_i \cap U_j$.
\item		Whenever $U_i \subset U_j$ for two uniformized sets, there is an injection
		$\lambda_{ij} : \{ V_i, G_i, \pi_i \} \rightarrow \{ V_j, G_j, \pi_j \}$.
\end{itemize}
If each chart in ${\mathcal F}$ is reduced, then $Q$ is said to be a {\bf reduced} orbifold.  Otherwise, $Q$ is {\bf unreduced}.
\end{definition}

The definition of an orbifold with boundary is identical, except that it allows
orbifold charts with boundary.  In this case,
$\partial Q := \{ \pi_i (\partial V_i) : \{ V_i , G_i, \pi_i \} \in {\mathcal F} \}$ is the
{\bf boundary} of $Q$.  It is easy to see that the restrictions
$\{ \partial V_i , G_i, (\pi_i)_{| \partial V_i} \}$ endow $\partial Q$ with the structure
of an orbifold.
We will sometimes refer to an orbifold as an {\bf orbifold without boundary} for emphasis.

It is easy to see that, given an unreduced orbifold $Q$, one can associate to it
a reduced orbifold $Q_{red}$ by redefining the group in each chart to be $G_i/K_i$,
where $K_i$ again denotes the kernel of the $G_i$-action on $V_i$.

Fix $p \in Q$, and say $p \in U$ for some set $U \subset Q$ uniformized by
$\{ V, G, \pi \}$.  Let $\tilde{p} \in V$ such that $\pi(\tilde{p}) = p$, and
let $I_{\tilde{p}}$ denote the isotropy subgroup of $\tilde{p} \in V$.
The isomorphism class of $I_{\tilde{p}}$ depends only on $p$; indeed, if
$\tilde{p}^\prime$ is another choice of a lift, then there is a group element
$\gamma \in G$ such that $\gamma \tilde{p}^\prime = \tilde{p}$, so that
$I_{\tilde{p}}$ and $I_{\tilde{p}^\prime}$ are conjugate via $\gamma$.
Similary, if $U^\prime, \{ V^\prime, G^\prime , \pi^\prime \}$ is another choice
of chart with $p \in U^\prime$ (and we assume, by shrinking domains if necessary, that
$U \subseteq U^\prime$), then there is an injection
$\lambda : \{ V, G, \pi \} \rightarrow \{ V^\prime, G^\prime, \pi^\prime \}$
with associated homomorphism $f : G \rightarrow G^\prime$ that maps $I_{\tilde{p}}$
isomorphically onto $I_{\phi (\tilde{p}) }$ (see \cite{satake2}, page 468).
We will often
refer to the (isomorphism class) of this group as the {\bf isotropy group of $p$},
denoted $I_p$.  If $I_p \neq 1$ (with respect to the orbifold structure of $Q_{red}$ in the case that $Q$ is not reduced),
then $p$ is {\bf singular}; otherwise, it is {\bf nonsingular}.  The collection of
singular points of $Q$ is denoted $\Sigma_Q$.

% XXXXXXXXXXXXXXXXXXXXXXXXXXXX	SUBSUBSECTION: EXAMPLES

\subsection{Examples}
\label{examples}

Before proceeding, we give some examples of orbifolds.

% XXXXXXXXXXXXXXXXXXXXXXXXXXXX	EXAMPLE

\begin{example}

If $M$ is a smooth manifold and $G$ a  group that acts properly discontinuously
on $M$ such that the fixed point set of each element of $G$ has codimension $\leq 2$,
then the quotient $M/G$ is an orbifold (\cite{satake2}, \cite{thurston}).

% XXXXXXXXXXXXXXXXXXXXXXXXXXXX	CLAIM

Similarly, we have the following claim.

\begin{claim}
\label{mwbmfg}
If $M$ is a manifold with boundary and $G$ a group with action
as above that fixes the boundary, then $M/G$ is an orbifold with boundary.
\end{claim}

This is proven as follows, following Thurston's proof (\cite{thurston})
for the case of manifolds without boundary.

Let $p \in Q$, and let $\tilde{p} \in M$ project to $p$.  Note that, by
hypothesis, if $\tilde{p} \in \partial M$, then any point projecting to
$p$ is also an element of $\partial M$.  Hence, we may define
$\partial Q$ to be the set of points which are projections of points in
$\partial M$.

Now, in the case that $p$ is not an element of $\partial Q$, we may use
Thurston's proof (by intersecting the neighborhood of $\tilde{p}$ with the interior
of $M$) to produce a chart  near $p$.  If $p \in \partial Q$, then let $I_{\tilde{p}}$
denote the isotropy group of $\tilde{p}$.  There is a neighborhood $V_{\tilde{p}}$
of $\tilde{p}$ invariant by $I_{\tilde{p}}$ and disjoint from its translates by elements
of $G$ not in $I_{\tilde{p}}$.  Note that $V_{\tilde{p}}$
is diffeomorphic to a neighborhood
of $\R_+^n$, where $n$ is the dimension of $M$, and that
$\partial \R_+^n$ is clearly $G$-invariant.  The
projection $U_p = V_{\tilde{p}} / I_{\tilde{p}}$ is by fiat a
homeomorphism.

To obtain a suitable cover of $M/G$, augment some cover $\{ U_p \}$ by
adjoining finite intersections.  Whenever
$U_{p_1} \cap \cdots \cap U_{p_k} \neq \emptyset$, this means some set of translates
$\gamma_1 V_{p_1} \cap \cdots \cap \gamma_k V_{p_k}$ has a
corresponding non-empty intesection.  This intersection may be taken to be
$U_{p_1} \cap \cdots \cap U_{p_k}$ with associated group
$\gamma_1 I_1 \gamma_1^{-1} \cap \cdots \cap \gamma_k I_k \gamma_k^{-1}$.

With this, we need only note that Thurston's proof (applied to the boundary,
and restricted to $\partial M$) will give a cover of $\partial Q$, consistent
with that on $Q$, making $\partial Q$ an orbifold.  This completes the proof.

\end{example}

Orbifolds which arise in this way are called {\bf global quotient orbifolds}
or {\bf global quotients}.  They are examples of {\bf good orbifolds,} orbifolds
whose universal cover is a manifold or, equivalently, orbifolds diffeomorphic
to the quotient of a manifold by the proper action of a discrete group.
Orbifolds which are not good are called {\bf bad orbifolds}. 

% XXXXXXXXXXXXXXXXXXXXXXXXXXXX	EXAMPLE

\begin{example}
(Kawasaki \cite{kawasaki1})
If $G$ is a Lie group that acts smoothly on a smooth manifold $M$ such that
\begin{itemize}
\item		$\forall x \in M$, the isotropy subgroup $G_x$ is compact,
\item		$\forall x \in M$, there is a smooth slice $S_x$ at $x$,
\item		$\forall x, y \in M$ such that $y \notin Gx$, there are slices
		$S_x$ and $S_y$ with $GS_x \cap GS_y = \emptyset$, and
\item		the dimension of $G_x$ is constant on $M$,
\end{itemize}
then $M/G$ is an orbifold.  Moreover, every $n$-dimensional orbifold $Q$ can be expressed
in this manner with $M$ the orthonormal frame bundle of $Q$ and $G = O(n)$.
\end{example}

% XXXXXXXXXXXXXXXXXXXXXXXXXXXX	EXAMPLE

\begin{example}
\label{expt}
The simplest example of an orbifold is that of a point $p$ with the trivial action
of a finite group $G$.  Note that these orbifolds are not reduced unless $G = 1$.
\end{example}

% XXXXXXXXXXXXXXXXXXXXXXXXXXXX	EXAMPLE

\begin{example}
\label{extear}
The $\Z_k$-teardrop is a well-known example of a bad orbifold, whose
underlying topological space is $S^2$.  It has one singular point, covered by
a chart $\{ V_1, G_1, \pi_1 \}$ where $V_1$ is an open disk in $\R^2$ about the
origin on which $G_1 = \Z_k$ acts via rotations, and every other point is covered
by a chart with the trivial group (see Figure~\ref{teardrop})).

% XXXXXXXXXXXXXXXXXXXXXXXXXXXX	FIGURE

\begin{figure}[h]
\centerline{
\includegraphics{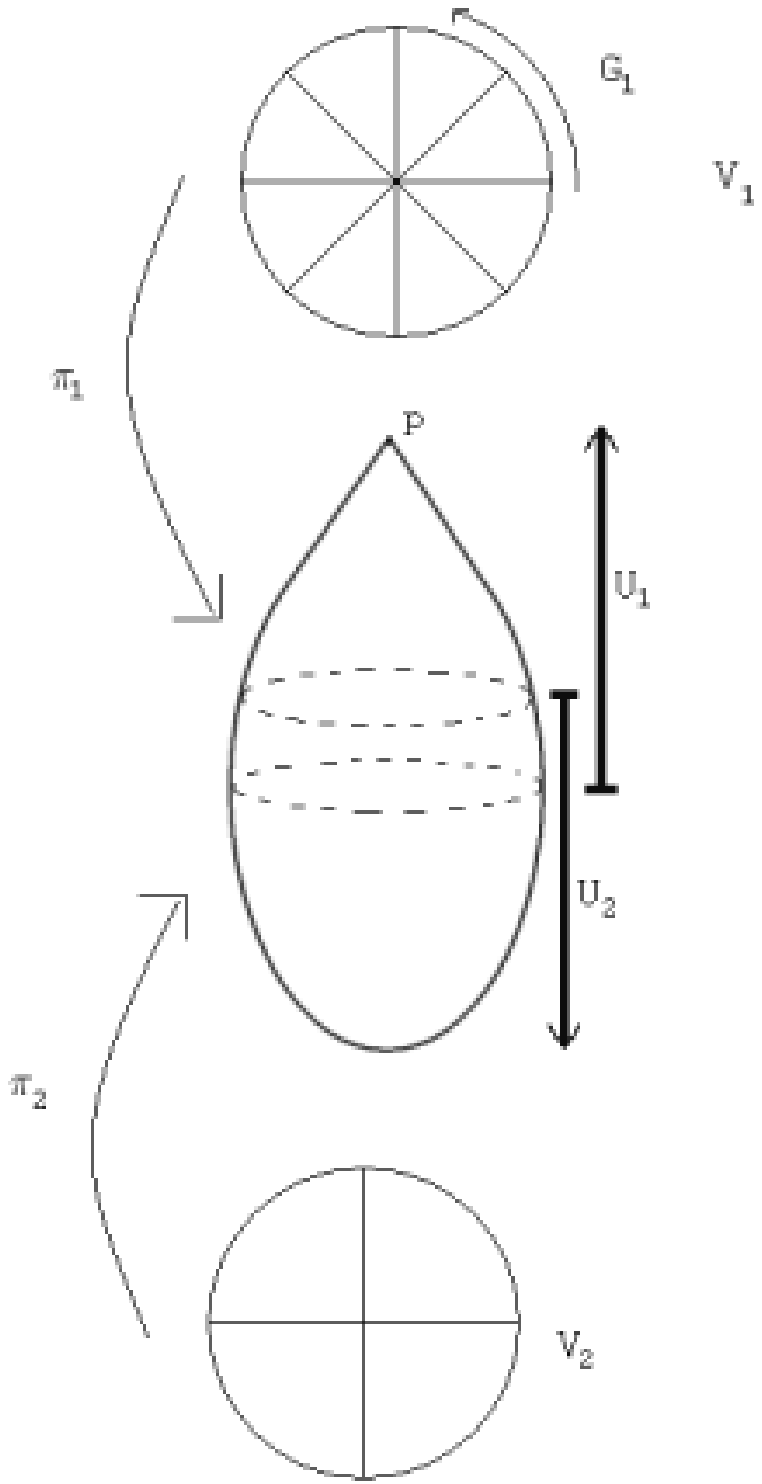}
}
\caption{\label{teardrop} 
{\footnotesize The $\Z_k$-teardrop.}}
\end{figure}

In fact, every compact 2-dimensional orbifold (with or without boundary) can be constructed from
a compact 2-dimensional manifold by removing disks and replacing them with $\R^2 /G$ for some
finite cyclic group $G < SO(2)$ with one singular point.  Note that in dimension 2,
the boundary of the orbifold cannot contain singular points.
\end{example}

The underlying space of an orbifold need not be a topological manifold, even in dimensions 3,
as is demonstrated by the following example (taken from \cite{stanhope}).

% XXXXXXXXXXXXXXXXXXXXXXXXXXXX	EXAMPLE

\begin{example}
\label{exrp2}
Let $\Z_2$ act on $\R^3$ via the antipodal map, and then the origin is the only fixed point.
Then $\R^3 / \Z_2$ is clearly an orbifold with one singular point, yet its underlying space
is homeomorphic to a cone on $\R \mathbb{P}^2$, which is not a manifold.
\end{example}

% XXXXXXXXXXXXXXXXXXXXXXXXXXXX	EXAMPLE

\begin{example}[The $\Z_k$-$\Z_l$-Solid Hollow Football]
\label{exsolidhollow}
An example of a bad orbifold with  boundary is the {\bf
$\Z_k$-$\Z_l$-solid hollow football} with $k \neq l$.  It is homeomorphic to the
manifold with boundary $\{ (x,y,z) \in \R^3 :  1 \leq \sqrt{x^2 + y^2 + z^2} \leq 2 \}$
with two boundary components, both homeomorphic to $S^2$.  Its orbifold structure,
however, is such that both of the boundary components have the orbifold structure
of the $\Z_k$-$\Z_l$-football (the sphere with two singular points having respective
groups $\Z_k$ and $\Z_l$; see \cite{borzellino}, page 16, Example 10),
and the interior has two singular sets, both homeomorphic to
a line segment, with isotropy subgroups of order
$k$ and $l$, respectively (see Figure~\ref{klfball}).

%	XXXXXXXXXXXXXXXXXXXXXX	Figure of Z_k-Z_l solid hollow football

\begin{figure}[h]
\centerline{
\includegraphics{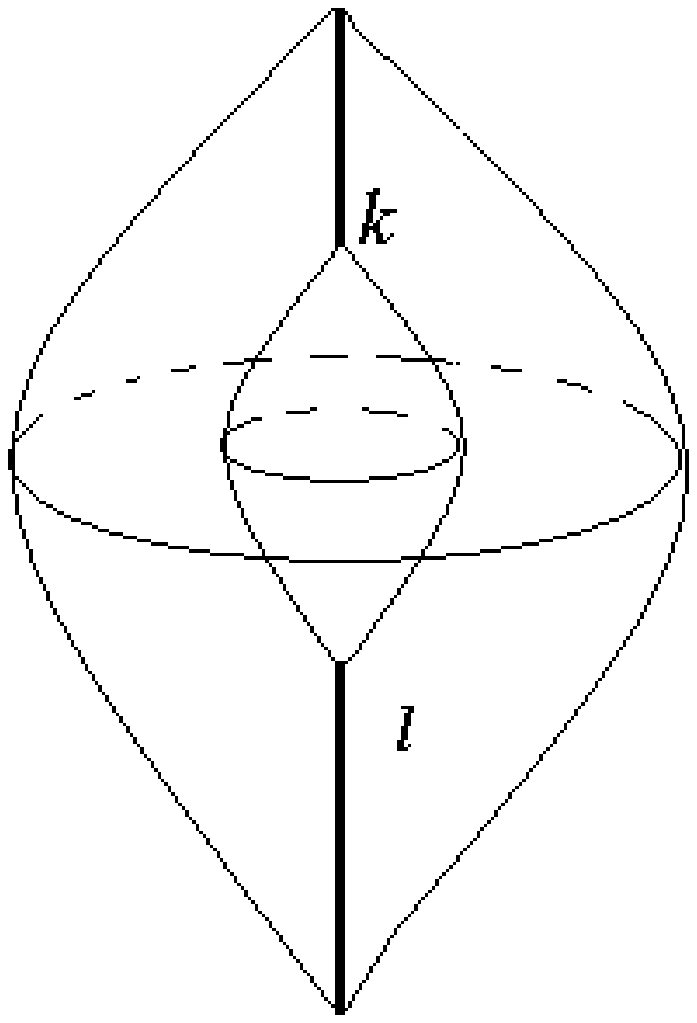}
}
\caption{\label{klfball} 
{\footnotesize The $\Z_k$-$\Z_l$-solid hollow football is the region
between the two closed 2-dimensional orbifolds.}
}
\end{figure}
\end{example}

Note that this orbifold is good when $k = l$.
For in this case, the orbifold can be expressed
as $M/\Z_k$ where $M = \{ (x,y,z) \in \R^3 :  1 \leq \sqrt{x^2 + y^2 + z^2} \leq 2 \}$
and $\Z_k$ acts via rotation about the $z$-axis.  
That this orbifold with boundary is bad whenever $k \neq l$ follows from the fact that
the $\Z_k$-$\Z_l$-football is bad and the following proposition.

% XXXXXXXXXXXXXXXXXXXXXXXXXXXX	PROPOSITION

\begin{proposition}

Let $Q$ be a good orbifold with  boundary.  The $\partial Q$ is a good orbifold.

\end{proposition}

\noindent {\it Proof:}

If $Q = M/G$ for some manifold $M$ and group $G$, then
$\partial Q = \partial M / (G_{\mid \partial M})$.

\bigskip \noindent Q.E.D. \bigskip

% XXXXXXXXXXXXXXXXXXXXXXXXXXXX	SUBSUBSECTION: STRUCTURES ON ORBIFOLDS

\subsection{Structures on Orbifolds}
\label{structures}

The next step is to introduce the appropriate notion of a vector bundle on an orbifold.
The following definition follows \cite{liang} (compare \cite{satake2} and
\cite{ruangwt}; note that our definition of an orbifold vector bundle corresponds
to Ruan's definition of a {\bf good} orbifold vector bundle).

% XXXXXXXXXXXXXXXXXXXXXXXXXXXX	DEFINITION

\begin{definition}[orbifold vector bundlle]
Let $Q$ be a connected orbifold.  By an {\bf orbifold vector bundle $E$ of rank $l$}, we mean
a collection consisting, for each set $U \subset Q$ uniformized by $\{ V, G, \pi \}$,
of a $G$-bundle $E_V$ over $V$ of rank $l$ such that the $G$-action
on $V$ and $E_V$ have the same kernel.  We require that for each injection
$\lambda_{12} : \{ V_1 , G_1 , \pi_1 \} \rightarrow \{ V_2 , G_2 , \pi_2 \}$,
there is a bundle map
$\phi_{12}^\ast : (E_{V_2})_{\mid \phi_{12}(V_1)} \rightarrow E_{V_1}$
such that if $\phi_{12} \circ \gamma_1 = \gamma_2 \circ \phi_{12}$ for some $\gamma_i \in G_i$,
then $\gamma_1 \circ \phi_{12}^\ast = \phi_{12}^\ast \circ \gamma_2$.
The {\bf total space} of the bundle $E$, also denoted $E$, is formed from the
collection $E_V /G$ by identifying points $(p, v) \in E_{V_1}$ and
$(q, w) \in E_{V_2}$ whenever there is an injection
$\lambda_{12} : \{ V_1 , G_1 , \pi_1 \} \rightarrow \{ V_2 , G_2 , \pi_2 \}$
such that $q = \phi_{12}(p)$ and $\phi_{12}^\ast (w) = v$.  It is clear that the total
space of a bundle $E$ is an orbifold.  Moreover, if $\rho_V$ denotes the
projection $\rho_V : E_V \rightarrow V$ for each $V$, then the collection
of these projections patch together to form a well-defined map
$\rho : E \rightarrow Q$, called the {\bf projection} of the bundle $E$.
In the case where $Q$ is not connected, we will require that the rank of the
$E_V$ is constant only on the connected components of $Q$.

By a {\bf section} $s$ of an orbifold bundle $E$, we mean a collection
$s_V$ of sections of the $G$-bundles $\{ E_V : \{ V, G, \pi \}$ is an orbifold chart
for $Q$ $\}$ such that if
$\lambda_{ij} : \{ V_i , G_i , \pi_i \} \rightarrow \{ V_j , G_j , \pi_j \}$
is an injection, then
$\phi_{ij}^\ast s_{V_j} [ \phi_{ij} (\tilde{p})] = s_{V_i} (\tilde{p})$ for each
$\tilde{p} \in V$.  We require that if $\tilde{p} \in V$ is fixed by
$\gamma \in G$, then $s(\tilde{p})$ is also fixed by $\gamma$.

It is clear that this collection defines a well-defined map $s : Q \rightarrow E$
from $Q$ into the total space of the bundle $E$ such that $\rho \circ s$ is the
identity on $Q$.
\end{definition}

It is important to notice that $E$ is not generally a vector bundle over $Q$, as
the fiber $\rho^{-1}(p)$ over a point $p \in Q$ is not always a vector space.
In general, $\rho^{-1}(p) \cong \R^l / I_p$, which is a vector space only
when $I_{\tilde{p}}$ acts trivially; i.e. when $p$ is nonsingular.

In particular, the {\bf tangent bundle} $\rho : TQ \rightarrow Q$ of an orbifold is
defined to be the collection of tangent bundles $TV$ over each $V$ where the $G$-action
is given by the Jacobian of the injection induced by the group action.  More
explicitly, each $\gamma \in G$ induces an injection
$\lambda_\gamma : \{ V, G, \pi \} \rightarrow \{ V, G, \pi \}$, as was explained in
Subsection \ref{orbandwbd}.
For $\tilde{p} \in V$ and $p \in U$ such that $\pi (\tilde{p}) = p$, let
\[
	g_\gamma (\tilde{p})
	=
	\left( \frac{ \partial (u^i \circ \gamma)}
			{\partial u^j} \right)
\]
be the Jacobian matrix of the action of $\gamma$ at $\tilde{p}$ for a fixed choice
$\{ u^i \}$ of coordinates for $V$.  Then 
$\{ g_\gamma (\tilde{p}): \gamma \in I_{\tilde{p}} \}$ defines an
$I_{\tilde{p}}$-action on $T_{\tilde{p}} V_i$ (see \cite{satake2}).
Note that the space
of vector fields on a uniformized set $U$, i.e. of sections of the orbifold
tangent bundle over $U$, corresponds
to the space of $G$-invariant vector fields on $V$. If $p \in Q$, we
refer to the maximal vector space in $T_pQ := \rho^{-1}(p)$ as the
{\bf space of tangent vectors at $p$.}  The tangent vectors at
a singular point are tangent to the singular set.

In particular, given a metric on $TQ$ (i.e. a $G$-invariant metric
on $Q$), the exponential map exp$:TQ \rightarrow Q$ is defined.  Using
this map and the fact that the $G$-action on the tangent bundle is linear, given
an orbifold chart $\{ V, G, \pi \}$, we define
an equivalent orbifold chart $\{ V^\prime, G^\prime, \pi^\prime \}$ such that
$U = U^\prime$ (i.e. $\pi(V) = \pi^\prime(V^\prime)$), and such that
$G^\prime$ acts on $V^\prime \subseteq \R^n$ as a subgroup of $O(n)$.  Picking a point
$p \in U^\prime$ and reducing domains, then, we always have that there is an
orbifold chart $\{ V^{\prime\prime}, G^{\prime\prime}, \pi^{\prime\prime} \}$
for an open neighborhood of $U^{\prime\prime}$ of $p$ such that
$G^{\prime\prime}$ acts linearly on $V^{\prime\prime}$ and
$p = \pi^{\prime\prime}(0)$.  In particular, $G^{\prime\prime}$ is the
isotropy subgroup of $p$.  Throughout, we will
use the convention that, for a given point $p \in Q$, a chart labeled
$\{ V_p, G_p, \pi_p \}$ has these properties with respect to $p$.  We refer to such
a chart as an {\bf orbifold chart at $p$ } (see \cite{ruangwt} for details).

In the same manner, we can define the cotangent bundle, its exterior powers, etc.
of an orbifold $Q$.  Differential forms are integrated over orbifolds in the following sense:
If $\omega$ is a differential form on $Q$ whose support is contained in
a uniformized set $U$, then the integral $\int\limits_Q \omega$ is defined to be
\[
	\frac{1}{|G|} \int\limits_V \tilde{\omega},
\]
where $\tilde{\omega}$ is the pullback of $\omega$ via the projection $\pi$.
Note that this integral does not depend on the choice of orbifold chart.
The integral of a globally-defined differential form is then defined in the
same way, using a partition of unity subordinate to a cover consisting of
uniformized sets (see \cite{satake2}).

% XXXXXXXXXXXXXXXXXXXXXXXXXXXXXXXXXXXXXXXXXXXXXXXXXXXXXXXXXXXXXXXXXXX
%			Section: Singular Dimension
% XXXXXXXXXXXXXXXXXXXXXXXXXXXXXXXXXXXXXXXXXXXXXXXXXXXXXXXXXXXXXXXXXXX

\section{The Dimension of a Singularity}

% XXXXXXXXXXXXXXXXXXXXXXXXXXXX	Subsection: Motivation

In this section, we are interested in understanding the behavior of vector fields on
orbifolds and orbifolds with boundary.  The primary restriction will be the dimension
changes in the space of tangent vectors, the maximal vector space contained in a fiber
of the orbifold tangent bundle.

\subsection{Motivation and Definition}

Throughout this section,
let $Q$ be an $n$-dimensional reduced orbifold, $p \in Q$, and let
$\{ V_p, G_p, \pi_p \}$ be an orbifold chart at $p$.  As was noted above, the tangent
bundle $TQ$ of $Q$ is not generally a vector bundle.  Indeed, if $p$ is a singular
point of $Q$ (i.e. $G_p \neq 1$), then letting $\rho : TQ \rightarrow Q$
again denote the projection,
$\rho^{-1} (p) \cong \R^n / \{ g_\gamma (\tilde{p}) : \gamma \in G \}$
for any lift $\tilde{p}$ of $p$ in $V_p$; the fiber is not a vector space
(recall that $g_\gamma (\tilde{p})$ denotes the infinitesimal action of
$\gamma \in G_p$ on $T_{\tilde{p}}V_p$; see Subsection \ref{structures}).
This fiber, however, is larger than the space of tangent vectors at $p$ (the vectors
in $T_{\tilde{p}} V_p$ which are fixed by $g_\gamma$ for each $\gamma \in G_p$),
i.e. the largest vector space contained in $\rho^{-1}(p)$.  For a
vector field $X$ on $Q$, $X(p)$ is always an element of the maximal vector space
contained in $T_p Q$.

Therefore, the space of tangent vectors of $Q$ at $p$ may, as a vector
space, have a smaller dimension than the dimension $n$ of the
orbifold.  This will play an important role in understanding the zeros
of vector fields on orbifolds.  In particular, if the space of tangent
vectors has dimension 0 at any point, then any vector field must
clearly vanish at that point; this is much different than the case of a manifold,
and motivates the following definition.

% XXXXXXXXXXXXXXXXXXXXXXXXXXXX	DEFINITION 

\begin{definition}[Dimension of a Singularity]

With the setup as above, we say that $p \in Q$ has {\bf singular dimension $k$}
(or that $p$ is a {\bf singularity of dimension $k$}) if the space of
tangent vectors at $p$ has dimension $k$.

\end{definition}

First, we verify that this definition is well-defined.

% XXXXXXXXXXXXXXXXXXXXXXXXXXXX	PROPOSITION

\begin{proposition} \label{dimsingcharindy}

The singular dimension of a point $p$ does not depend on the choice
of the orbifold chart, nor on the choice of the lift $\tilde{p}$ of $p$
in the chart.

\end{proposition}

\noindent	{\it Proof:}

Fix $p \in Q$, and let $\{ V_i, \pi_i, G_i \}$ and
$\{ V_j, \pi_j, G_j \}$ be orbifold charts (with $U_i := \pi_i (V_i)$ and
$U_j := \pi_j (V_j)$ as usual) such that $p \in U_i \cap U_j$.
Suppose first that $U_i \subseteq U_j$, and then by the definition of an
orbifold, there is an injection $\lambda_{ij}$ with embedding $\phi_{ij}: V_i \rightarrow V_j$.

Now, let $\tilde{p}_i$ be a point in $V_i$ such that
$\pi_i(\tilde{p}_i) = p$, and then $\tilde{p}_j := \phi_{ij}(\tilde{p}_i)$
has the property that $\pi_j(\tilde{p}_j) = p$ (by the definition of
$\phi_{ij}$).  With respect to the chart $V_i$, the
singular dimension of $p$ is the dimension of the space of tangent vectors at
$p$; i.e. the dimension of the vector subspace of $T_{\tilde{p}_i} V_i$
which is fixed by the infinitessimal $I_{\tilde{p}_i}$-action.  We have that
$\phi_{ij}$ is a diffeomorphism of $V_i$ onto an open subset of $V_j$.
Moreover, as the groups $I_{\tilde{p}_i}$ and $I_{\tilde{p}_j}$ are
isomorphic, we have that the associated injective homomorphism
\[
	f_{ij} : G_i \rightarrow G_j
\]
maps $I_{\tilde{p}_i}$ onto $I_{\tilde{p}_j}$.
Hence, for any $v \in T_{\tilde{p}_i} V_i$ such that
$g_\gamma (\tilde{p}_i) v = v$ for every $\gamma \in I_{\tilde{p}_i}$,
it is clear that $g_{\gamma}(\tilde{p}_i) d(\phi_{ij})_{\tilde{p}_i}(v) =
d(\phi_{ij})_{\tilde{p}_i}(v)$ for every $\gamma \in I_{\tilde{p}_j}$,
and conversely.  Therefore, the map $d(\phi_{ij})_{\tilde{p}_i}$ restricts
to a vector space isomorphism of the fixed-point set of the $I_{\tilde{p}_i}$-action
on $T_{\tilde{p}_i} V_i$ onto the fixed-point set of the $I_{\tilde{p}_j}$-action on
$T_{\tilde{p}_j} V_j$, ensuring that their dimensions are equal.

That the singular dimension at $p$ does not depend on our choice of the lift
$\tilde{p}$ of $p$ is now clear: if $\tilde{p}_i^\prime$ is another choice of
a lift in $V_i$, then there is an injection $\lambda_{ii} = (\phi_{ii}, f_{ii})$ of
$(V_i, G_i, \pi_i)$ into itself such that $\phi_{ii}(\tilde{p}_i) = \tilde{p}_i^\prime$.
Hence we may apply the above argument.

In the case where $U_i$ is not a subset of $U_j$, there is a chart over some set
$U_k$ with $p \in U_k \subseteq U_i \cap U_j$.
The above argument gives us that the dimension is the same with respect to
$U_k$ as with respect to $U_i$, and the same with respect to
$U_k$ as with respect to $U_j$.

\bigskip \noindent Q.E.D. \bigskip

% XXXXXXXXXXXXXXXXXXXXXXXXXXXX	Subsection: Properties

\subsection{Properties}

Note that every point in $Q$ has a singular dimension, not simply the
singular points.  In particular, we have

% XXXXXXXXXXXXXXXXXXXXXXXXXXXX	PROPOSITION

\begin{proposition}	\label{nonsingisdimn}
The non-singular points in $Q$ are precisely the points with
singular dimension $n$, where $n$ is the dimension of $Q$.
\end{proposition}

\noindent	{\it Proof:}

Let $p$ be a point in $Q$, and suppose that
$p$ has singular dimension $n$.  Let $\{ V_p, G_p, \pi_p \}$ be
an orbifold chart at $p$.  Taking a lift
$\tilde{p} \in V_p$ such that $\pi(\tilde{p}) = p$, we then
have by hypothesis that the space of $I_{\tilde{p}}$-invariant
vectors in $T_{\tilde{p}} V_p$, has dimension $n$, and hence that
$I_{\tilde{p}}$ acts trivially on $T_{\tilde{p}} V_p$.  Of course, as $Q$ is
reduced, and as $G_p = I_{\tilde{p}}$ by our choice of charts,
this implies that $G_p = 1$, and that $p$ is not a singular point.

Conversely, for $p$ a non-singular point in $Q$, $I_{\tilde{p}} = 1$
for any lifting $\tilde{p}$ of $p$ into any chart, and hence fixes the entire
$n$-dimensional vector space $T_{\tilde{p}} V$ in the corresponding chart.

\bigskip \noindent Q.E.D. \bigskip

Recall that $\Sigma_Q$ denotes the set of all singular points
in $Q$ (i.e. points $p \in Q$ such that $G_p \neq 1$ for any chart
$\{ V_p, G_p, \pi_p \}$ at $p$).
For each $k$ with $k=0, 1,  \ldots , n$, let $\Sigma_k$ denote
the set of all singular points with singular dimension $k$.
Clearly $\Sigma_Q = \bigcup\limits_{k=0}^{n-1} \Sigma_k$.  As was pointed out in
Satake \cite{satake2}, $Q \backslash \Sigma_Q$ (which, in our notation
is $\Sigma_n$) is an $n$-dimensional manifold.  However, we also have
following proposition (see also \cite{kawasaki1}):

% XXXXXXXXXXXXXXXXXXXXXXXXXXXX	PROPOSITION

\begin{proposition}	\label{singkmanifold}
Let $Q$ be an $n$-dimensional orbifold (with boundary).  Then for
$k = 0, 1, \ldots , n$, $\Sigma_k$ naturally has the structure of a
$k$-dimensional manifold (with boundary, and
$\partial \Sigma_k = \Sigma_k \cap \partial Q$).
\end{proposition}

\noindent	{\it Proof:}

Fix $k$ with $0 \le k \le n$, let $p$ be a point in $\Sigma_k$, and fix
an orbifold chart $\{ V_p, G_p, \pi_p \}$ at $p$.
Let $\tilde{p} \in V_p$ be a lift of $p$ in the chart.  By the
definition of $\Sigma_k$, we have that the $I_{\tilde{p}}$-invariant subspace
of $T_{\tilde{p}} V_p$ is of dimension $k$.  Hence, as $G_p$ acts
linearly on $V$, there is a $k$-dimensional subspace $Y$ of
$\R^n \supset V_p$ on which $I_{\tilde{p}}$ acts trivially.

If $k = 0$, then $Y = \{ \tilde{p} \}$ is clearly the only invariant point of
$V_p$, and hence $U_p$ is an open set containing $p$ in which
$p$ is the only 0-dimensional singularity.  Therefore, the $0$-dimensional
singularities are isolated and form a 0-manifold.  The rest of the proof
will deal with the case $k > 0$.

First suppose that $p$ is not a boundary point of $Q$ (so, in particular,
we may assume that $V_p$ does not have boundary).  Note that by the definition
of $Y$, each point $\tilde{x} \in Y$ is fixed by
$I_{\tilde{p}}$.  Hence, the isotropy subgroup $I_{\tilde{x}}$ of each
$\tilde{x} \in Y$ in $G_p$ clearly contains $I_{\tilde{p}}$.  However,
by our choice of charts, $I_{\tilde{p}} = G_p$, so that each
$I_{\tilde{x}} = G_p$.  Hence, we may take $W = Y \cap V_p$ to be an open neighborhood
of $\tilde{p}$ in $Y$ with constant isotropy.

Note that as $W \subset Y$, the restriction $(\pi_p)_{\mid W}$ of
$\pi_p$ to $W$ is an injective map of $W$ onto $\pi_p (W)$.
Hence, as $\pi_p(W) \cap \Sigma_k$ is clearly open in $\Sigma_k$,
$(\pi_p)_{\mid W}$ is a ${\mathcal C}^\infty$
bijection of $W$ onto a neighborhood of $p$ in $\Sigma_k$.
In particular, as $W$ is open in $Y = \R^k$, $\pi_p$
restricts to a manifold chart for $\Sigma_k$ near $p$.
That such charts cover $\Sigma_k$ is clear, as $p$ was arbitrary.
Moreover, that such charts intersect with the appropriate transition
follows directly from the existence of the injections $\lambda_{ij}$
and associated $\phi_{ij}$.

Similarly, if $Q$ has boundary and $p \in \partial Q$, then $Y$ is a subspace of
$\R^n \supset \R_+^n \supset V_p$, and $W = Y \cap V_p$ is an open subset
of $Y \cap \R_+^n$.  Hence, we see that $\pi_p$ restricts to a chart
of a manifold with boundary on $\Sigma_k$ with $p \in \partial \Sigma_k$.

\bigskip \noindent Q.E.D. \bigskip

% XXXXXXXXXXXXXXXXXXXXXXXXXXXX	EXAMPLE

\begin{example}
\label{exfig8}
We note here that the connected components of the set $\Sigma_Q$
need not have the structure of a manifold.
For example, within a 3-dimensional orbifold, we may have a singular set in the
shape of an 8 that fails to be a manifold at one point $p$.  This can happen,
for instance, when a chart $\{ V_1, G_1 , \pi_1 \}$ at $p$ has
group
\[
	G_1	=	\langle R_\pi^x , R_{\frac{2\pi}{3}}^z \rangle
\]
where
\[
	R_\pi^x = \left[ \begin{array}{ccc}
		1	&	0	&	0	\\\\
		0	&	-1	&	0	\\\\
		0	&	0	&	-1
			\end{array} \right]
\]
denotes a rotation of $\pi$ radians about the $x$-axis, and
\[
	R_{\frac{2\pi}{3}}^z = \left[ \begin{array}{ccc}
		-\frac{1}{2}		&	-\frac{\sqrt{3}}{2}	&	0	\\\\
		\frac{\sqrt{3}}{2}	&	-\frac{1}{2}		&	0	\\\\
		0	&	0	&	1
			\end{array} \right]
	=
	\left[ \begin{array}{ccc}
		\cos\left(\frac{2\pi}{3}\right)			&
		-\sin\left(\frac{2\pi}{3}\right)			&	0	\\\\
		\sin\left(\frac{2\pi}{3}\right)			&
		\cos\left(\frac{2\pi}{3}\right)			&	0	\\\\
		0	&	0	&	1
			\end{array} \right]
\]
denotes a rotation of $\frac{2\pi}{3}$ radians about the $z$-axis.
By direct computation,
$R_\pi^x R_{\frac{2\pi}{3}}^z (R_\pi^x)^{-1} = (R_{\frac{2\pi}{3}}^z)^{-1}$ so
that $G_1$ is isomorphic to the dihedral group $D_6$.  Define further the charts
$\{ V_2, G_2 , \pi_2 \}$ and $\{ V_3, G_3 , \pi_3 \}$ covering the rest of $\Sigma_Q$
with groups
\[
	G_2 = \langle \R_\pi^x \rangle \cong \Z_2
\]
and
\[
	G_3 = \langle \R_{\frac{2\pi}{3}}^z \rangle \cong \Z_3
\]
as pictured (see Figure~\ref{eight}).

% XXXXXXXXXXXXXXXXXXXXXXXXXXXX	FIGURE

\begin{figure}[h]
\centerline{
\includegraphics{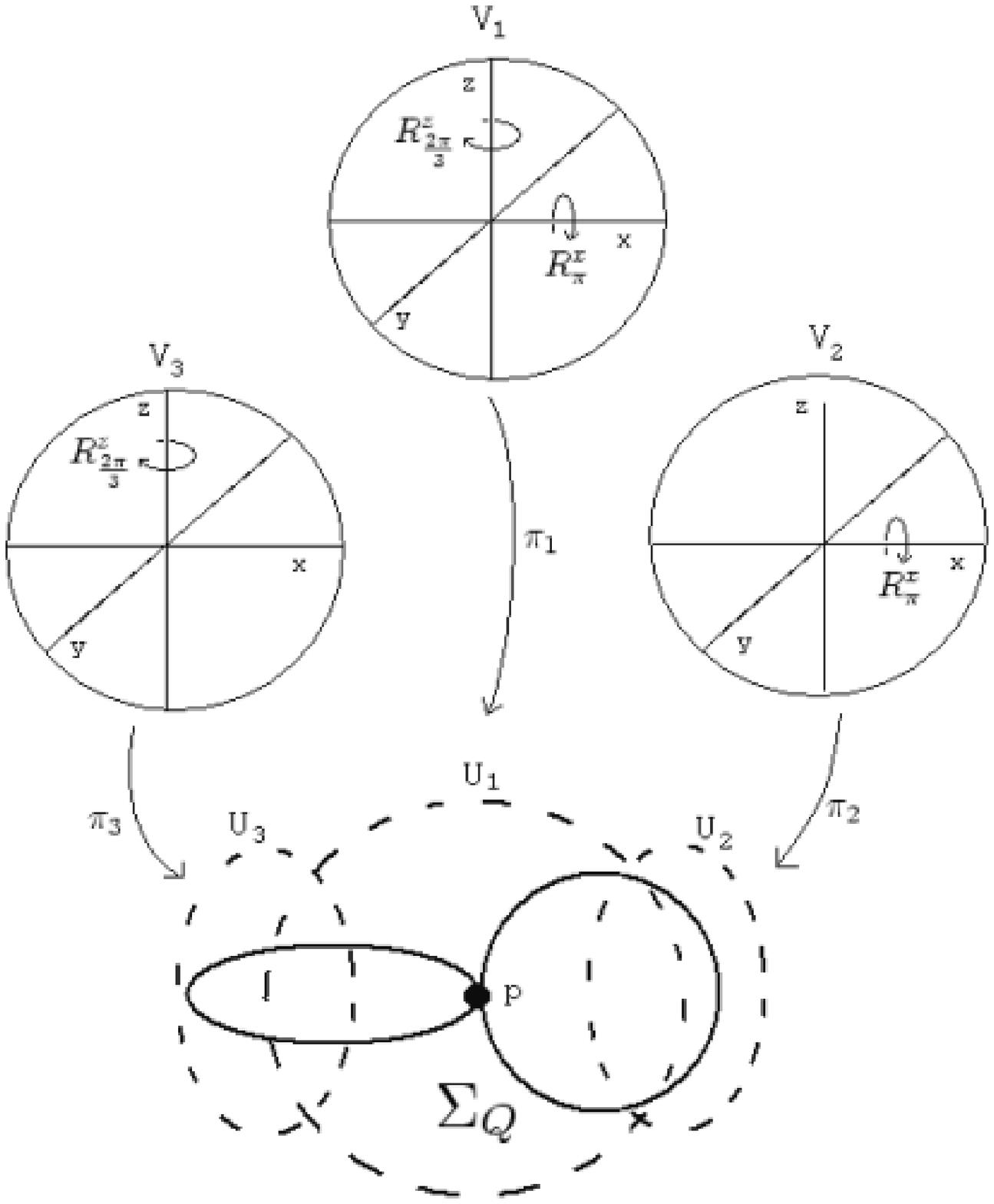}
}
\caption{\label{eight} 
{\footnotesize A singular set $\Sigma_Q$ in a 3-dimensional orbifold which is not
a manifold.}}
\end{figure}

We note, however, that the set $\Sigma_1 = \Sigma_Q \backslash \{ p \}$ is an
open 1-manifold diffeomorphic to $(0,1) \sqcup (2,3)$,
and the set $\Sigma_0 = \{ p \}$ is a 0-manifold.
\end{example}

This example serves to illutrate the structure of the set $\Sigma_Q$
almost in general.  For, intuitively, in any chart
$\{ V, G , \pi \}$, vectors tangent to $V$ which are
invariant under the $\{ g_\gamma : \gamma \in G \}$-action must clearly
be tangent to the
singular set $\Sigma_Q$.  Hence, at points where two singular strata `run into'
one another, the dimension of the space of tangent vectors decreases.
The intersection of the closures of two connected components of $\Sigma_k$,
for some $k$, belong to $\Sigma_j$ for some $j < k$.

% XXXXXXXXXXXXXXXXXXXXXXXXXXXX	COROLLARY

We state the following obvious corollaries to Proposition \ref{singkmanifold}
and its proof.

\begin{corollary}	\label{sing0finite}
With the notation as above, $\Sigma_0$ is a set of isolated points.  In
particular, if $Q$ is compact, $\Sigma_0$ is finite.
\end{corollary}

% XXXXXXXXXXXXXXXXXXXXXXXXXXXX	COROLLARY

\begin{corollary}	\label{tangentsing}
Let $p \in \Sigma_k \subset Q$.  Then the set of tangent vectors of $Q$ at $p$
is canonically identified with $T_p \Sigma_k$.
\end{corollary}

% XXXXXXXXXXXXXXXXXXXXXXXXXXXX	Subsection: Examples

\subsection{Examples of 0-Dimensional Singularities}

Our primary interest here is the case of singularities with dimension
zero, over which the largest vector space contained in $T_p Q$
is $\{ \mathbf{0} \}$.  For over such points, any vector field must vanish.

% XXXXXXXXXXXXXXXXXXXXXXXXXXXX	EXAMPLE

\begin{example}
Consider, for example, the cone $C$, which is the quotient of $\R^2$
by the group of rotations
\[
	G	=	\langle \gamma \rangle	=	\Z_3
\]
where
\[
	\gamma :=	\left[\begin{array}{cc}
				\cos \frac{2\pi}{3}	&	-\sin \frac{2\pi}{3}	\\\\
				\sin \frac{2\pi}{3}	&	\cos \frac{2\pi}{3}
			\end{array}\right] 
\]
(see Figure~\ref{Z3cone}).

% XXXXXXXXXXXXXXXXXXXXXXXXXXXX	FIGURE

\begin{figure}[h]
\centerline{
\includegraphics{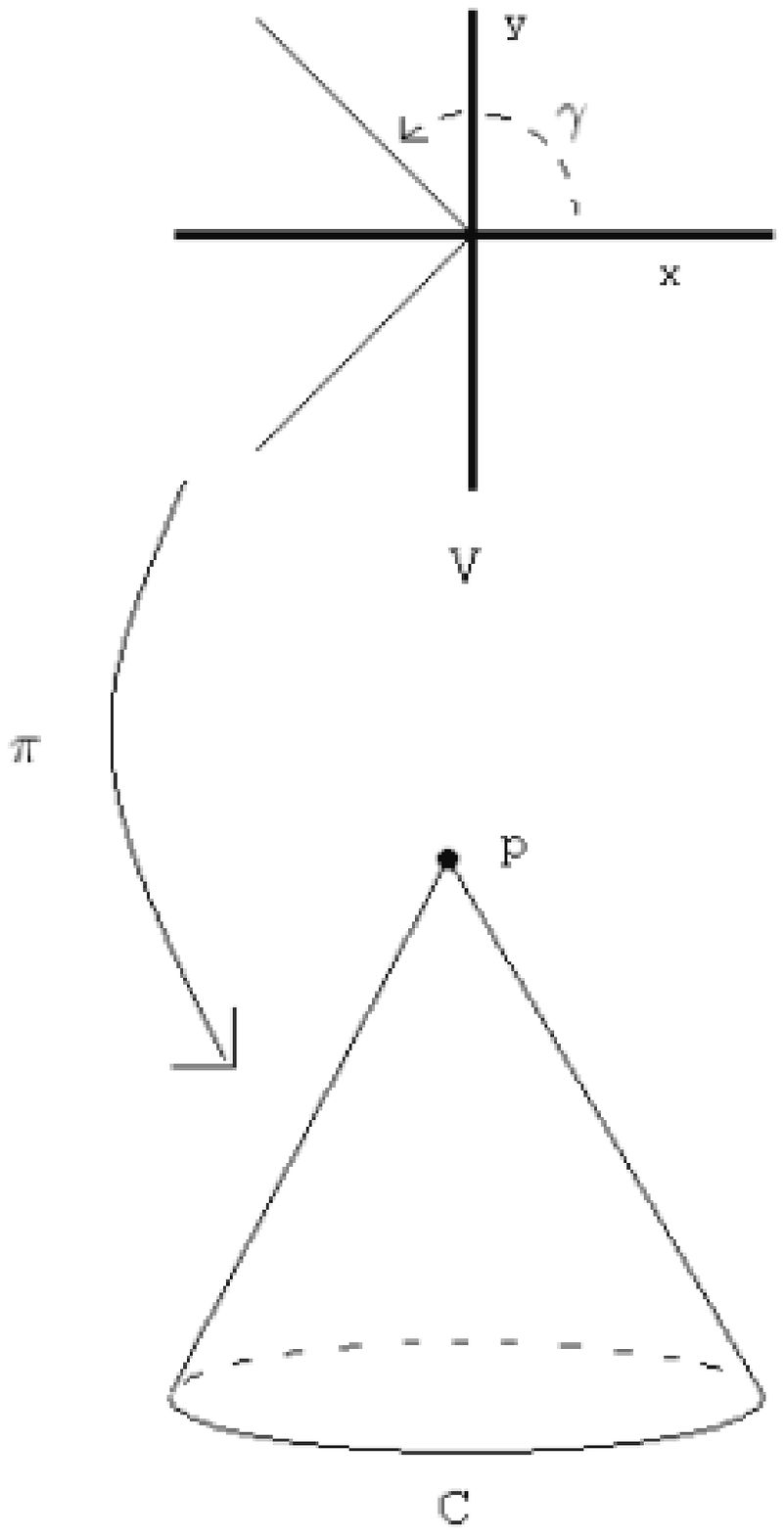}
}
\caption{\label{Z3cone} 
{\footnotesize The cone $C = \R^2 / \Z_3$.}}
\end{figure}

This is clearly an orbifold with one chart where $V = \R^2$, $G \cong \Z_3$, and
$\pi$ is the quotient map induced by the group action; and
one singular point $p = \pi(0,0)$.  A vector field on $C$ is precisely
a $G$-invariant vector field on $V = \R^2$.  However, as $\gamma$
has no nontrivial eigenvectors, any vector on the tangent space at $(0,0)$ is
not fixed by $G$, so that any vector field on $C$ must vanish at the
singular point.
\end{example}

It may be helpful if we develop the tangent
space of this previous example explicitly.  Any point $q \in C$, $q \neq p$,
with lift $\tilde{q} \in V$ has trivial isotropy
group $I_{\tilde{q}} = 1$.  Hence, the action of $I_{\tilde{q}}$ on the
tangent space $T_{\tilde{q}} V \cong \{ \tilde{q} \} \times \R^2$
is trivial, and the space of tangent vectors is
\[
\begin{array}{rcl}
	T_q Q
	&=&
	\pi^\ast [(T_{\tilde{q}}V)^{I_{\tilde{q}}}]				\\\\
	&=&
	\pi^\ast (T_{\tilde{q}}V)							\\\\
	&\cong&
	\R^2.
\end{array}
\]

The set $\Sigma_2 = C \backslash \{ p \}$ is clearly a smooth 2-dimensional manifold.
Moreover, taking a point $q \in \Sigma_2$ and an open ball $V^\prime$ about
$q$ which does not intersect $\tilde{p} = (0,0)$ (small enough so that
$\forall \gamma \in G, (\gamma V^\prime) \cap V^\prime = \emptyset$),
the restriction of $\pi$ to $V^\prime$ gives a manifold chart of $\Sigma_2$ near $q$,
and the manifold tangent space of $\Sigma_2$ at $q$ is exacly the
orbifold tangent space of $C$ at $q$.

Now, the point $\tilde{p} = (0,0)$ has isotropy group $I_{\tilde{p}} = G$.
Identifying $T_{\tilde{p}} V$ with $\R^2$ in the usual way,
we have that as $G$ is linear, and $g_\gamma(\tilde{p}) = \tilde{p}$
(identifying coordinates on $V$ with coordinates on $T_{\tilde{p}} V$ via the
exponential map;
this simply states that the Jacobian of a linear operator is itself).  Then
$(T_{\tilde{p}}V)^{I_{\tilde{p}}}$ is defined to be the set of vectors in
$T_{\tilde{p}} V \cong \R^2$ which are invariant under the action of
the group generated by
$g_\gamma(\tilde{p}) = g = \begin{array}{c} \\ \left[\begin{array}{cc}
		\cos \frac{2\pi}{3}	&	-\sin \frac{2\pi}{3}	\\\\
		\sin \frac{2\pi}{3}	&	\cos \frac{2\pi}{3}
	\end{array}\right] \\\\ \end{array}$,
which is clearly only the zero vector.  Hence,
$T_pQ = \pi^\ast (T_{\tilde{p}}V^{I_{\tilde{p}}}) = \pi^\ast( \{ \mathbf{0} \} ) = \{ \mathbf{0} \}$,
and the space of tangent vectors to $C$ at $p$ contains only the zero vector.

The following example is meant to illustrate the limitations of Satake's
definition of an orbifold with boundary (see \cite{satake2} for the definition).
In particular, it is an example in which there is no non-vanishing vector field
on the boundary.

% XXXXXXXXXXXXXXXXXXXXXXXXXXXX	EXAMPLE

\begin{example}
\label{slicedcone}
Now consider the orbifold $C_0$ with boundary (using Satake's definition of
an orbifold with boundary) defined to be the union of the following subsets of $C$:
\[
\begin{array}{rcl}
	V_1:	&=&
		\left\{ (\rho,\theta) : \rho \geq 0 ,
			\frac{-\pi}{6} \leq \theta \leq \frac{\pi}{6} \right\},	\\\\
	V_2:	&=&
		\left\{ (\rho,\theta) : \rho \geq 0 ,
			\frac{3\pi}{6} \leq \theta \leq \frac{5\pi}{6} \right\},
			\;\mbox{and}						\\\\
	V_3:	&=&
		\left\{ (\rho,\theta) : \rho \geq 0 ,
			\frac{7\pi}{6} \leq \theta \leq \frac{9\pi}{6} \right\} 
\end{array}
\]
(we use standard polar coordinates $(\rho, \theta)$ on $\R^2$; see
Figure~\ref{Z3conewithbd}).

% XXXXXXXXXXXXXXXXXXXXXXXXXXXX	FIGURE

\begin{figure}[h]
\centerline{
\includegraphics{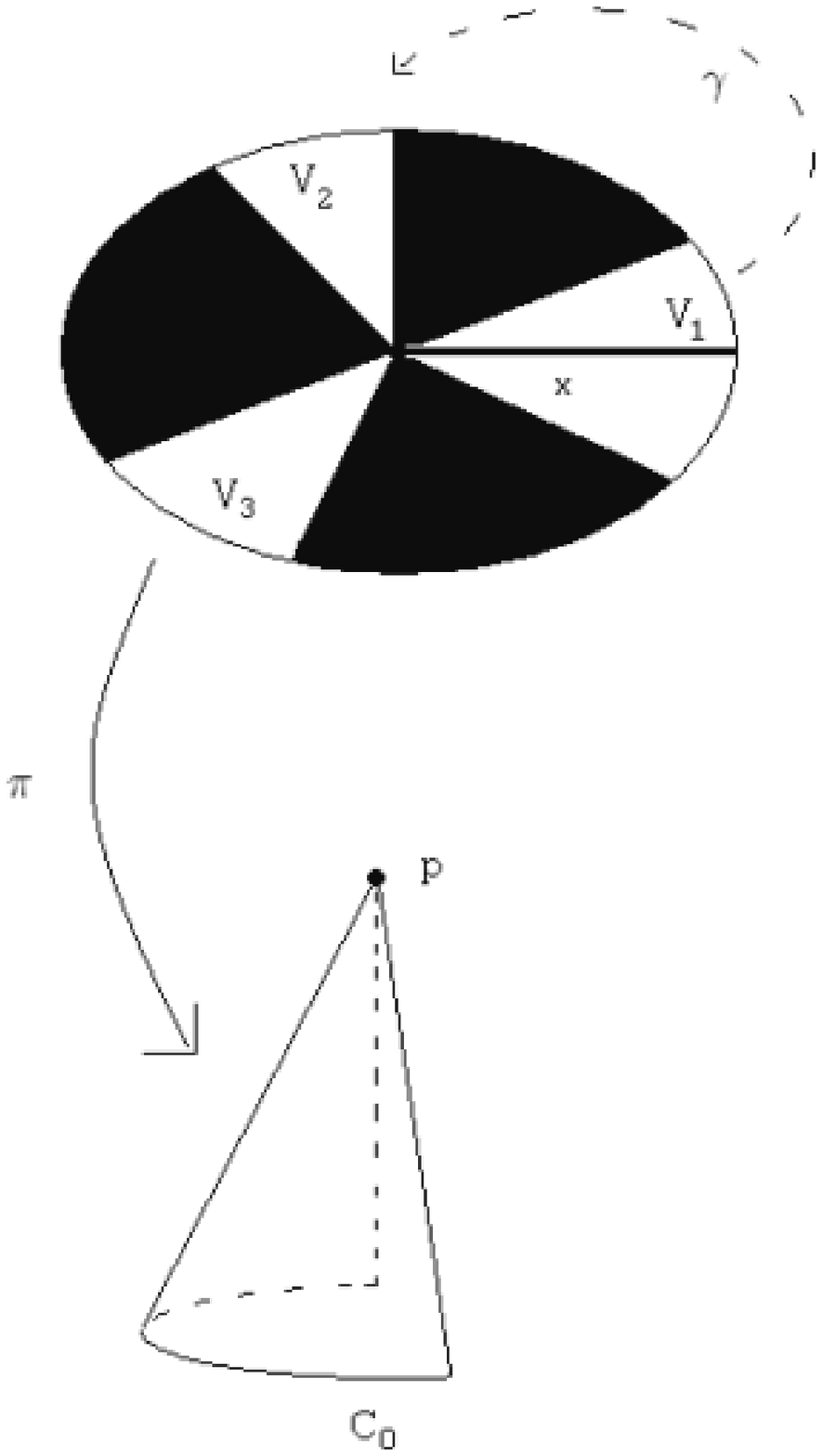}
}
\caption{\label{Z3conewithbd} 
{\footnotesize A dimension-0 singularity on the boundary of an orbifold (using Satake's original
definition of an orbifold with boundary).}}
\end{figure}

Then
\[
	C_0: = \pi(V_1) \cup \pi(V_2) \cup
		\pi(V_3)
\]
is clearly an orbifold with boundary containing one singular point
$p = \pi(0,0)$ of dimension zero.  Hence, as $p \in \partial C_0$,
there is no vector field $X$ on $Q_0$ which is non-vanishing on the
boundary.  Indeed, for any such vector field, $X(p) = \mathbf{0}$.
Note further that $C_0$ does not admit a product structure near the boundary;
i.e. there is no open neighborhood of the boundary in $Q$ diffeomorphic to
$[0, \epsilon) \times \partial Q$.
\end{example}

We should note that these two previous examples are non-compact for simplicity
of exposition, but that, by reducing the domain of the charts and patching them
together with charts where $G$ is the trivial group, the
same type of singularity can clearly occur in the case of a compact orbifold:
the first as the sole singularity in the $\Z_3$-teardrop, and the
second in a teardrop which is missing a piece homeomorphic to a disc,
where the singularity occurs on the boundary.

For our definition of orbifolds, however, the above cannot occur.  Indeed,
it is trivial to show (using a chart with linear group action) that
if $Q$ is an orbifold with boundary $M$, then a neighborhood of $M$ in $Q$
is diffeomorphic to $[0, \epsilon) \times M)$.  Therefore,
for every point $p$ on the boundary of $Q$,
the tangent space $T_p Q$ contains a one-dimensional subspace independent from
$T_p M$.  In particular, although $M$ may, as an orbifold, contain 0-dimensional
singularities, as a subset of $Q$, it does not.
As well, by Corollary \ref{sing0finite}, on a compact orbifold, we need not
worry about the existence of a vector field with a finite number of zeros.

%% file: thesisch3.tex
% Thesis Chapter 3, Christopher Seaton, seatonc@colorado.edu
% XXXXXXXXXXXXXXXXXXXXXXXXXXXXXXXXXXXXXXXXXXXXXXXXXXXXXXXXXXXXXXXXXXX
% XXXXXXXXXXXXXXXXXXXXXXXXXXXXXXXXXXXXXXXXXXXXXXXXXXXXXXXXXXXXXXXXXXX
%			Chapter 3: The Gauss-Bonnet Theorem
% XXXXXXXXXXXXXXXXXXXXXXXXXXXXXXXXXXXXXXXXXXXXXXXXXXXXXXXXXXXXXXXXXXX
% XXXXXXXXXXXXXXXXXXXXXXXXXXXXXXXXXXXXXXXXXXXXXXXXXXXXXXXXXXXXXXXXXXX
%
%	section 1:	Introduction for this chapter
%	section 2:	Theorems on the level of forms
%									setup
%									notes on the integrand definitions
%									results of satake for closed orbifolds
%									calculations for outer unit normal vector field
%									proofs of the theorems
%	section 3:	Thom Isomorphism for Orbifolds
%									thom iso in de Rham cohomology
%									the case of a global quotient
%									the case of a general orbifold
%	section 4:	Theorems on the level of cohomolgy
%

\chapter{The First Gauss-Bonnet and Poincar\'{e}-Hopf Theorems for
	Orbifolds With Boundary}
\label{chgbph}

% XXXXXXXXXXXXXXXXXXXXXXXXXXXXXXXXXXXXXXXXXXXXXXXXXXXXXXXXXXXXXXXXXXX
%			Section: Introduction
% XXXXXXXXXXXXXXXXXXXXXXXXXXXXXXXXXXXXXXXXXXXXXXXXXXXXXXXXXXXXXXXXXXX

\section{Introduction}

Our goal in this section is to generalize Satake's Gauss-Bonnet Theorem for
orbifolds with boundary \cite{satake2} using the more modern definition of
an orbifold with boundary.  We will see that, with this definition,
Satake's boundary term can be simplified considerably (indeed, in some
cases it will vanish).

To this end, a note is necessary.  The statemenet of Satake's Gauss-Bonnet
Theorem for orbifolds with boundary refers to an outward-pointing unit
normal vector field on the boundary $M$ of $Q$.  In the case that $Q$
is taken to be an orbifold with boundary as defined by Satake, the boundary may very
well contain 0-dimensional singular points (as in the case of the sliced
cone in Example \ref{slicedcone}),
in which case this vector field does not exist---indeed, in such
a case there will be no non-vanishing vector field on the boundary.  However,
our definition does not allow such behavior.

This section follows constructions in \cite{chern1}, \cite{satake2}, and \cite{sha}.
However, our notation primarily follows that of \cite{sha}.  Where possible, we
state our results and constructions for general orbifold vector bundles, though
our primary application will be to the tangent bundle.  In particular, many of the
definitions will be made in as a general a context as possible to facilitate the proof
of the Thom Isomorphism Theorem for Orbifolds (Theorem \ref{thomgeneral}) in
the general case.  All orbifolds under consideration are assumed to be
reduced.

% XXXXXXXXXXXXXXXXXXXXXXXXXXXXXXXXXXXXXXXXXXXXXXXXXXXXXXXXXXXXXXXXXXX
%			Section: Results in terms of Differential Forms
% XXXXXXXXXXXXXXXXXXXXXXXXXXXXXXXXXXXXXXXXXXXXXXXXXXXXXXXXXXXXXXXXXXX

\section{The Results on the Level of Differential Forms}

% XXXXXXXXXXXXXXXXXXXXXXXXXXXX	SUBSUBSECTION: SETUP

\subsection{The Setup}
\label{setup}

In this subsection, we fix the context and notation for the rest of the section.
Let $Q$ be a compact, connected, orientable orbifold of dimension $n$ with boundary
$\partial Q =: M$, taking the orientation of $M$ to be the orientation 
inherited from
$Q$ with respect to the outer unit normal vector field.
We will use with liberty the fact that that there is a
neighborhood of $M$ in $Q$ diffeomorphic to $[0, \epsilon) \times M$.
Let $\rho : E \rightarrow Q$ be an orbifold vector bundle over $Q$ of rank
$l = 2m$ or $2m + 1$ equipped with a Euclidean metric.  Let $SE$  denote
the unit sphere bundle of $E$ with respect to the metric, with projection
still denoted $\rho$. Fix a compatible $SO(n)$-connection $\omega$ with
curvature $\Omega$ on $E$.  With respect to a local orthonormal frame
$(e_1 , \ldots , e_l)$ for $E$, we let $(u_1 , \ldots , u_l)$ denote the
component functions on $E$ and $(\theta_1 , \ldots , \theta_l )$ the
basic forms: $\theta:= du + \omega u$.

For the specific case where $E = TQ$ is the tangent bundle of $Q$,
we require that the metric respect the product structure of $Q$ near the boundary.
In particular, at any point $p \in M$, with respect to any chart
$\{ V, G, \pi \}$ and any oriented orthonormal frame field $(e_1, e_2, \ldots, e_n)$
for the fiber $(\pi^\ast TQ)_{\tilde{p}}$ at $\tilde{p}$ ($\pi(\tilde{p}) = p$ as usual),
we have that $\omega_{i,j}= \Omega_{i,j} = 0$ whenever $i = 1$ or $j = 1$.  Note that
in the case of the tangent bundle, $l = n$.

Let $X$ be a vector field on $Q$ which is non-zero on $M$ and has a finite number of
fixed points $p_1, \ldots , p_s$ on the interior of $Q$.
Denote by $B_r(p_i)$ the geodesic ball about $p_i$ of radius $r$, and for ease of
notation, let $B_r(P) := \bigcup\limits_{i=1}^s B_r(p_i)$.
Let $\alpha : Q \backslash \{p_1, \ldots , p_s \} \rightarrow ST$ denote
the section of the sphere bundle $STQ$ induced by $X$.
We will sometimes require that $X$ extend the outer unit normal vector field on $M$,
in which case we will use the notation $X_0$ for the vector field and $\alpha_0$ for
the induced section of the sphere bundle.

% XXXXXXXXXXXXXXXXXXXXXXXXXXXX	SUBSUSECTION: INTEGRAND NOTES

\subsection{Notes on the Definitions of the Integrands}

A note is in order on our definition of the Euler curvature form $E(\Omega)$
and its secondary form $\Psi$.  These forms, which were originally defined
by Chern \cite{chern1} in the case of the tangent bundle of a closed Riemannian
manifold of
even dimension, were used by Satake \cite{satake2} and Sha \cite{sha} in
generalizations of the Gauss-Bonnet Theorem to the case of closed orbifolds and of the
Poincar\'{e}-Hopf Theorem to the case of manifolds with boundary, respectively.
However, the definitions of the two forms differ slightly, primarily due to Chern's
further developments in \cite{chern2}.  In particular, in \cite{chern2}, Chern extended
the definition of these forms to odd-dimensional manifolds and reversed the sign of the
Euler curvature form.  Hence, some care must be taken with respect to how the various
definitions of the forms fit together.

In \cite{chern1}, the forms first appeared on the tangent bundle of a closed
Riemannian manifold of dimension $n = 2m$.  They were called
$\Omega$ and $\Pi$, respectively, and defined as follows with respect to a frame
$(e_1, \ldots , e_n)$ (again with component functions
$(u_1, \ldots , u_n)$ and basic forms $(\theta_1, \ldots , \theta_n)$):
\[
	\Omega := (-1)^{m-1} \frac{1}{2^{2m} \pi^m m! }
		\sum\limits_{\tau \in S(n)} (-1)^\tau
		\Omega_{\tau(1)\tau(2)} \wedge \cdots \wedge
		\Omega_{\tau(n-1)\tau(n)}
\]
is the Euler curvature form, and
\[
	\Pi := \frac{1}{\pi^m} \sum\limits_{k = 0}^{m - 1}
		\frac{1}{1 \cdot 3 \cdots (n - 2k - 1)} \Phi_k
\]
its secondary form on the unit sphere bundle, where
\[
	\Phi_k = \sum\limits_{\tau \in S(n)} (-1)^\tau u_{\tau(1)}
		\theta_{\tau(2)} \wedge \cdots \wedge
		\theta_{\tau(n - 2k)} \wedge
		\Omega_{\tau(n-2k+1)\tau(n-2k+2)}
		\wedge \cdots \wedge \Omega_{\tau(n - 1)\tau(n)}
\]
for $k = 0, 1, \ldots , m - 1$ (as usual, $S(n)$ denotes the group of
permutations of $n$ letters).

In \cite{chern2}, Chern introduced a new definition
of the form $\Pi$ on both even and odd-dimensional manifolds:
\[
	\Pi := 
	\left\{
	\begin{array}{ll}
		\frac{1}{\pi^m} \sum\limits_{k=0}^{m-1} (-1)^k
		\frac{1}{1 \cdot 3 \cdots (n - 2k - 1)2^{m + k}
		\cdot k !} \Phi_k , 	&\mbox{if} \; n = 2m
			\; \mbox{is even,}					\\\\
		\frac{1}{2^n \pi^m m!}
		\sum\limits_{k=0}^{m} (-1)^{k + 1}
		{ m \choose k }
		\Phi_k				&\mbox{if} \; n = 2m + 1
		\; \mbox{is odd,}
	\end{array}
	\right.
\]
where
\[
	\Phi_k = \sum\limits_{\tau \in S(n)} \Omega_{\tau(1) \tau(2)}
		\wedge \cdots \wedge \Omega_{\tau(2k - 1) \tau (2k)}
		\wedge \omega_{\tau(2k+1) n} \wedge \cdots \wedge
		\omega_{\tau(n-1) n}
\]
for $k = 0, 1, \ldots , \left[ \frac{n}{2} \right] - 1$.  It is
pointed out (\cite{chern2}, page 675) that when $n$ is even, the
$\Phi_k$, and hence $\Pi$, reduce to their definitions in \cite{chern1}.
As was mentioned above, the definition of $\Omega$
differs from the previous only by a minus sign (so that a sign is introduced
in the relationship between the two forms).

The reason for our difficulty is that various authors have since used different
conventions regarding the signs of these forms,
so that the definitions of the forms are by no means standard.
The definitions of Chern in \cite{chern2} are such that, in the odd dimensional
case, the integral of the secondary form on a fiber of the sphere bundle is
$-1$ in the odd-dimensional case, and hence relates to the negative index of the
vector field.  This, of course, offers no difficulty for closed manifolds, as
the index of the vector field is in this case always zero.  However, in the case
of manifolds with boundary, this leads to a negative sign in the formula, which
is undesirable.
The reader who may compare our computations to those in \cite{chern1}, \cite{chern2},
\cite{satake2}, \cite{schwerdtfeger}, or \cite{sha} is {\bf warned} to compare the
definitions carefully and take into account the appropriate sign conventions.

We will follow Sha's notation for the most part, letting $E(\Omega)$ denote
the Euler curvature form (with respect to the connection $\Omega$) and $\Psi$ its
secondary form.  Our sign conventions are chosen so that in the resulting formula,
the index of the vector field has the same sign in both the even and odd cases. 
Hence, we will use the following definitions, for an arbitrary vector bundle of
rank $l$ (again, $l = n$ in the case of the tangent bundle):
\[
	E(\Omega) :=
	\left\{
	\begin{array}{ll}
		\frac{1}{2^{2m} \pi^m m! }
		\sum\limits_{\tau \in S(l)} (-1)^\tau
		\Omega_{\tau(1)\tau(2)} \wedge \cdots \wedge
		\Omega_{\tau(l-1)\tau(l)},
			&\mbox{if} \; l = 2m
			\; \mbox{is even,}					\\\\
		0,				&\mbox{if} \; l = 2m + 1
			\; \mbox{is odd}
	\end{array}
	\right.
\]
is the Euler curvature form, which agrees with the definition of Sha \cite{sha}.
The secondary form on the unit sphere bundle is
\[
	\Psi = 
	\left\{
	\begin{array}{ll}
		\frac{(-1)^m}{\pi^m} \sum\limits_{k=0}^{m-1} (-1)^k
		\frac{1}{1 \cdot 3 \cdots (l - 2k - 1)2^{m + k}
		\cdot k !} \Phi_k , 	&\mbox{if} \; l = 2m
			\; \mbox{is even,}					\\\\
		\frac{-1}{2^l \pi^m m!}
		\sum\limits_{k=0}^{m}
		{ m \choose k}
		\Phi_k ,				&\mbox{if} \; l = 2m + 1
			\; \mbox{is odd,}
	\end{array}
	\right.
\]
where
\[
	\Phi_k = \sum\limits_{\tau \in S(l)} (-1)^\tau u_{\tau(1)}
		\theta_{\tau(2)} \wedge \cdots \wedge
		\theta_{\tau(l - 2k)} \wedge
		\Omega_{\tau(l-2k+1)\tau(l-2k+2)}
		\wedge \cdots \wedge \Omega_{\tau(l - 1)\tau(l)} .
\]
We have that on the unit sphere bundle, $d\Psi = -\rho^\ast E(\Omega)$, where
$\rho : Q \rightarrow SE$ again denotes the bundle projection.  Moreover, we have
that $\int\limits_{SE_p} \Psi = \frac{1}{|I_p|}$ where $p$ is any point in $Q$ with
isotropy $I_p$.
Of course, these relations are preserved up to the sign using the definitions in
any of our references (in the case of a manifold, $\frac{1}{|I_p|}$ is
always equal to $1$).

% XXXXXXXXXXXXXXXXXXXXXXXXXXXX	SUBSUBSECTION: SATAKE RESULTS

\subsection{The Result of Satake for Orbifolds With Boundary}

In \cite{satake2}, we have the following Gauss-Bonnet Theorem.

\begin{theorem}(Satake)
\label{gb}
Let $Q$ be an oriented compact Riemannian orbifold with boundary $M$, let
${\mathcal N}$ be the outward-pointing unit normal vector field on $M$,
and let $\alpha_0$ be the induced section of $STQ_{|M}$.  Then we have
\[
	\int\limits_Q E(\Omega) = \chi_{orb}^\prime (Q) -
		\int\limits_M \alpha_0^\ast \Psi ,
\]
where $M$ is oriented by the induced orientation with respect to the outer normal
vector field on $M$.
\end{theorem}

It is noted that the requirement of orientability can be lifted by
using the standard reduction to the orientable double-cover.

Here, $\chi_{orb}(M)$ and $\chi_{orb}^\prime (Q)$ are the {\bf orbifold Euler
characteristic} and {\bf inner orbifold Euler characteristic} of $M$ and $Q$,
respectively, whose definitions we recall (see Appendix \ref{appa} for a discussion
of alternate Euler characteristics for orbifolds).

Satake's definition of the Euler characteristic for a closed orbifold without
boundary came from his proof of the Gauss-Bonnet Theorem in this case,
with the equation
\[
	\int\limits_Q E(\Omega) = \sum\limits_{i=1}^s \mbox{ind}_{X}(p_i).
\]
Since the left-hand side does not depend on the vector field $X$, and
the right-hand side does not depend on the metric, this number is
an invariant of the orbifold itself.  We denote this invariant
$\chi_{orb} (Q)$ and refer to it as the {\bf orbifold Euler
Characteristic}.  Satake goes on to show how this number can be computed
in terms of a suitable triangulation, the existence of such a triangulation
since having been demonstrated in \cite{moerdijk}.  Specifically, if ${\mathcal T}$
is a triangulation of $Q$ such that the order of the isotropy group is a
constant function on the interior of any simplex $\sigma \in {\mathcal T}$,
then letting $N_\sigma$ denote this order, we have
\[
	\chi_{orb} (Q) = \sum\limits_{\sigma \in {\mathcal T}}
		\frac{(-1)^{\mbox{dim}\:\sigma}}{N_\sigma}.
\]

Now, in the case that $Q$ has boundary, the {\bf inner orbifold Euler
Characteristic } $\chi_{orb}^\prime(Q)$ is defined similarly, but with
respect to a particularly chosen vector field: one which extends the
outward unit normal ${\mathcal N}$ on the boundary $M$ of $Q$.  In
this case, given a simplicial decomposition as above, we have
\[
	\chi_{orb}^\prime (Q) = \sum\limits_{\sigma \in {\mathcal T}_0}
		\frac{(-1)^{\mbox{dim}\:\sigma}}{N_\sigma},
\]
where ${\mathcal T}_0$ denotes the collection of simplices which are not
completely contained in the boundary.

% XXXXXXXXXXXXXXXXXXXXXXXXXXXX	SUBSUBSECTION: COMPUTATIONS FOR OUTER VECTOR FIELDS

\subsection{Computations of $\alpha^\ast(\Psi)$ in the Case That $X$ Extends
	the Outer Unit Normal Vector Field on $M$}

As was pointed out in Sha \cite{sha} for the case of manifolds, in the specific case
that the vector field $X_0$ extends the outward-pointing normal vector field
${\mathcal N}$ and $n$ is even, if $\alpha_0$ denotes the section of $ST Q$
induced by $X_0$, we have that
\[
	\alpha_0^\ast (\Psi)		=			0
\]
on $M$.  This is proven as follows:

With a chart $\{ V, G, \pi \}$
and an orthonormal oriented frame field $(e_1, e_2, \ldots, e_n)$ chosen
for a point $\tilde{p}$ on the boundary $\pi^{-1}(M) = \partial V$ such that the
$(e_2, \ldots , e_n)$ are tangent to the boundary of $V$,
the form $\Psi$ is a sum of terms of the form
\[
	\Phi_k = \sum\limits_{\tau \in S(n)} (-1)^\tau u_{\tau(1)}
		\theta_{\tau(2)} \wedge \cdots \wedge
		\theta_{\tau(n - 2k)} \wedge
		\Omega_{\tau(n-2k+1)\tau(n-2k+2)}
		\wedge \cdots \wedge \Omega_{\tau(n - 1)\tau(n)}
\]
for $k = 0, 1, \ldots m - 1$.
We have that the ${\mathcal N}$ is locally equal to $-e_1$ in the lift to $V$; it has
coordinates $(u_1, u_2, \ldots , u_n) = (-1,0,\ldots,0)$ with respect to this frame.
Hence, as $\theta = du + \omega u$, we have that $\alpha_0^\ast (\theta_i)$
vanishes for each $i > 1$.  Note that each term either has a factor of $\theta_i$
for $i > 1$, or has a factor of $\theta_1$ and a factor of $u_i$ for $i > 1$,
in which case $u_i = 0$.  Therefore, each of these terms vanish when composed with
$\alpha_0$, and hence
\[
	\alpha_0^\ast (\Psi) = 0.
\]

It will be worthwhile to see how these computations simplify
in the case of $n = 2m + 1$ odd (this is also noted by Sha in \cite{sha}).
As above, all of the $\theta_i$ factors vanish except for $\theta_1$,
but there does exist one $\Psi_k$ which does not
contain any such factors:
\[
\begin{array}{rcl}
	\Phi_m
		&=&
		\sum\limits_{\tau \in S(n)} (-1)^\tau u_{\tau(1)}
		\Omega_{\tau(n-2m+1)\tau(n-2m+2)}
		\wedge \cdots \wedge \Omega_{\tau(n - 1)\tau(n)}	\\\\
		&=&
		\sum\limits_{\tau \in S(n)} (-1)^\tau u_{\tau(1)}
		\Omega_{\tau(2)\tau(3)}
		\wedge \cdots \wedge \Omega_{\tau(n - 1)\tau(n)}.
\end{array}
\]

Moreover, in $\alpha_0^\ast (\Phi_m )$, we have that the coefficient
$u_{\tau(1)} = 0$ in every term except those such that
$\tau(1) = 1$ (recall that $u_1 = -1$ and $u_i = 0$ for $i > 1$).
Hence,
\[
\begin{array}{rcl}
	\alpha_0^\ast (\Phi_m )
		&=&
		u_1 \sum\limits_{\tau \in S(n - 1)} (-1)^\tau
		\Omega_{\tau(2)\tau(3)}
		\wedge \cdots \wedge \Omega_{\tau(n - 1)\tau(n)}	\\\\
		&=&
		-\sum\limits_{\tau \in S(n - 1)} (-1)^\tau 
		\Omega_{\tau(2)\tau(3)}
		\wedge \cdots \wedge \Omega_{\tau(n - 1)\tau(n)},
\end{array}
\]
where $S(n-1)$ is understood to be the group of permutations on
$\{ 2, 3, \ldots, n \}$.
So in this case,
\[
\begin{array}{rcl}
	\alpha_0^\ast (\Psi )
	&=&
	\frac{-1}{2^n \pi^m m!}
	{ m \choose m}
		\alpha_0^\ast (\Phi_m)						\\\\
	&=&
		\frac{-1}{2^n \pi^m m!}
		\left[ -\sum\limits_{\tau \in S(n - 1)} (-1)^\tau 
		\Omega_{\tau(2)\tau(3)}
		\wedge \cdots \wedge \Omega_{\tau(n - 1)\tau(n)} \right]
										\\\\
	&=&
	\frac{1}{2}\left[ \frac{1}{2^{2m} \pi^m m!}
		\sum\limits_{\tau \in S(n - 1)} (-1)^\tau 
		\Omega_{\tau(2)\tau(3)}
		\wedge \cdots \wedge \Omega_{\tau(n - 1)\tau(n)} \right]
\end{array}
\]
Recall that $(e_2, \ldots , e_n)$ is a frame
field for $TM$, and that $M$ is of dimension $n - 1 = 2m$.  Therefore,
the above form is precisely $\frac{1}{2}$ times the Euler curvature
form for $TM$.  In summary,
\[
	\alpha_0^\ast (\Psi ) = \frac{1}{2}E(\Omega_{|M}) .
\]

Hence, in the odd case,
\[
\begin{array}{rcl}
	\int\limits_M \alpha_0^\ast (\Psi)
	&=&
	\int\limits_M \frac{1}{2} E(\Omega_M)		\\\\
	&=&
	\frac{1}{2} \chi_{orb}(M),
\end{array}
\]
the last equality following from Satake's Gauss-Bonnet Theorem for
closed orbifolds.

% XXXXXXXXXXXXXXXXXXXXXXXXXXXX	SUBSUBSECTION: Theorems With Boundary

\subsection{The Theorems for Orbifolds With Boundary}

With this, we may restate Satake's result for orbifolds with 
boundary.

\begin{theorem}[The First Gauss-Bonnet Theorem for Orbifolds with Boundary]
\label{gbb1}
Let $Q$ be a compact orbifold of dimension $n$ with boundary $M$, and let
$E(\Omega)$ be defined as above in terms of the curvature $\Omega$
of a connection $\omega$.  Then
\[
	\int\limits_Q E(\Omega)	=
	\left\{	\begin{array}{ll}
		\chi_{orb}^\prime(Q),
			&			n = 2m,				\\\\
		\chi_{orb}^\prime(Q) - \frac{1}{2} \chi_{orb} (M),
			&			n = 2m + 1.
	\end{array} \right.
\]
\end{theorem}

Note that as $E(\Omega)$ is defined to be zero in the case that $n$ is odd, we have the
familiar relation
\[
	\chi_{orb}^\prime (Q) = \frac{1}{2} \chi_{orb} (M).
\]

\begin{example}
For example, the solid $\Z_3$-football (i.e. the closed ball in $\R^3$ with the
usual action of $\Z_3$ via rotations; see \cite{borzellino})
is an orbifold with  boundary whose inner orbifold Euler Characteristic
is $\frac{1}{3}$.  Its boundary has orbifold Euler
Characteristic $\frac{2}{3}$.  This can be easily verified with a simplicial
decomposition, or using the fact that the the space is diffeomorphic to
$\overline{\mathbb{D}_3}/ \Z_3$ and the boundary $S^2/ \Z_3$.
\end{example}

We are now in the position to extend the Poincar\'{e}-Hopf Theorem
to the case of a compact oriented orbifold with  boundary.
Begin with the setup given in Section \ref{setup} with $E = TQ$ and $X$
a vector field with a finite number of zeros $p_1, \ldots , p_s$ on the
interior of $Q$.  We require only that $X$ does not vanish on the boundary.
The index of $X$ at $p_i$ is defined in a manner analogous to the integral;
if $\{ V_{p_i}, G_{p_i}, \pi_{p_i} \}$ is a chart at $p_i$ and $\tilde{X}$
denotes the lift of $X$ to $V_{p_i}$, then the index of $X$ at $p$ is
\[
	\mbox{ind}_{X}(p_i) :=
	\frac{1}{|G_{p_i}|} \mbox{ind}_{\tilde{X}} (\tilde{p_i}),
\]
(see \cite{satake2}).  Of course, $\frac{1}{|G_{p_i}|} = \frac{1}{|I_{p_i}|}$
for this choice of chart; i.e. $G_{p_i}$ is the isotropy group of $p_i$.

In both the even and odd cases, we have that
\[
	d\Psi = -\rho^\ast E(\Omega)
\]
on the unit sphere bundle of the tangent bundle.  Moreover, at each singular point
$p_i$, we have that
\[
	\lim\limits_{r \to 0^+} \int\limits_{\partial B_r(p_i)} \alpha^\ast \Psi
	=
	-\mbox{ind}_X (p_i).
\]
This follows from the fact that the integral $\int\limits_{ST_p Q} \Psi$ of
$\Psi$ over any fiber $ST_p Q$ of the unit sphere bundle is $\frac{1}{|I_p|}$.
The minus sign is due to the fact that the orientation that $\partial B_r(p_i)$
inherits as a component of the boundary of $Q \backslash B_r (P)$ is the opposite
orientation of that used in definition of the index (see for example \cite{gp} or
\cite{milnor2}).
Recall that $B_r(P)$ denotes the union $\bigcup\limits_{i=1}^s B_r(p_i),$ and set
$\mbox{ind}(X) := \mbox{ind}_X (P) := \sum\limits_{i=1}^s \: \mbox{ind}_X (p_i)$.

With this, based on Sha's proof of the Poincar\'{e}-Hopf Theorem for
manifolds with boundary, we have the following:

% XXXXXXXXXXXXXXXXXXXXXXXXXXXX	even case

If $n = 2m$ is even, then
\[
\begin{array}{rcl}
	\chi_{orb}^\prime (Q)
	&=&
	\int\limits_Q E(\Omega)																									\\\\
	&&\mbox{(by Theorem \ref{gbb1})}
																																		\\\\
	&=&
	\lim\limits_{r \to 0^+} \int\limits_{Q \backslash B_r(P)}
		\alpha^\ast \rho^\ast (E(\Omega))																\\\\
	&&\mbox{(as} \; \rho \alpha \; \mbox{is the identity map on $Q$)}
																																		\\\\
	&=&
	-\lim\limits_{r \to 0^+} \int\limits_{Q \backslash B_r(P)}
		d\alpha^\ast (\Psi)																							\\\\
	&&\mbox{(as $\rho^\ast E(\Omega) = -d\Psi$ on $STQ$)}							\\\\
	&=&
	- \lim\limits_{r \to 0^+}
		\int\limits_{\partial B_r(P)} \alpha^\ast (\Psi)
		- \int\limits_M \alpha^\ast(\Psi)																			\\\\
	&&\mbox{(by Stokes' Theorem)}																			\\\\
	&=&
	\mbox{ind}(X) - \int\limits_M \alpha^\ast(\Psi),
\end{array}
\]
so that
\[
	\mbox{ind}(X)
	=
	\chi_{orb}^\prime (Q) + \int\limits_M \alpha^\ast(\Psi).
\]

% XXXXXXXXXXXXXXXXXXXXXXXXXXXX	odd case

Similarly, if $n = 2m + 1$ is odd, then as $E(\Omega) = 0$,
\[
\begin{array}{rcl}
	0
	&=&
	\lim\limits_{r \to 0^+} \int\limits_{Q \backslash B_r(P)}
		\alpha^\ast \rho^\ast (E(\Omega))				\\\\
	&=&
	-\lim\limits_{r \to 0^+} \int\limits_{Q \backslash B_r(P)}
		d\alpha^\ast (\Psi)											\\\\
	&=&
	- \lim\limits_{r \to 0^+} \int\limits_{\partial B_r(P)} \alpha^\ast (\Psi)
	- \int\limits_M \alpha^\ast(\Psi)								\\\\
	&=&
	\mbox{ind}(X) - \int\limits_M \alpha^\ast(\Psi),
\end{array}
\]
so that
\[
	\mbox{ind}(X)
	=
	\int\limits_M \alpha^\ast(\Psi) .
\]

% XXXXXXXXXXXXXXXXXXXXXXXXXXXX	PROPOSITION

In summary, we state the following.

\begin{proposition}
\label{phdependent}
Let $Q$ be a compact orbifold of dimension $n$ with boundary $M$.  Let
$X$ be vector field on $Q$ which has a finite number of singularities,
all of which occuring on the interior of $Q$.  Then
\[
	\mbox{ind}(X)		=
	\left\{	\begin{array}{ll}
		\chi_{orb}^\prime (Q) + \int\limits_M \alpha^\ast(\Psi),
			&			n = 2m,				\\\\
		\int\limits_M \alpha^\ast(\Psi),
			&			n = 2m + 1.
	\end{array} \right.
\]
\end{proposition}

% XXXXXXXXXXXXXXXXXXXXXXXXXXXXXXXXXXXXXXXXXXXXXXXXXXXXXXXXXXXXXXXXXXX
%			Section: Thom Iso for Orbifolds
% XXXXXXXXXXXXXXXXXXXXXXXXXXXXXXXXXXXXXXXXXXXXXXXXXXXXXXXXXXXXXXXXXXX

\section{The Thom Isomorphism Theorem for Orbifolds}

In Section \ref{thrmcohom} below, we will show that, as is the case with manifolds,
the cohomology class of the form $\Psi$ is an invariant of $STQ_{\mid M}$, and hence
does not depend on the various choices made.
In order to characterize the cohomology class of $\Psi$ in $H^n(STQ_{| M})$,
we need to determine its relationship with the Euler and Thom classes of the tangent
bundle $TQ_{|M}$.  To this end, we develop the Thom Isomorphism Theorem
for orbifolds in de Rham cohomology (Theorem \ref{thomgeneral}).

% XXXXXXXXXXXXXXXXXXXXXXXXXXXX	SUBSUBSECTION: thom iso in de Rham cohom

\subsection{The Thom Isomorphism in de Rham Cohomology}

We begin by stating the following Theorem (taken from \cite{milnor}).
In this case, $E_0$ denotes the set of nonzero vectors in the vector bundle
$\xi$ with total space $E$ (so that
$E_0 := E \backslash \{ \:\mbox{zero section} \: \}$), and $F$
is a typical fiber (with $F_0 := F \cap E_0$, the set of nonzero vectors in $F$).

% XXXXXXXXXXXXXXXXXXXXXXXXXXXX	THEOREM

\begin{theorem}[Thom Isomorphism Theorem]
\label{milnorthom}

Let $\xi$ be an oriented $l$-plane bundle with total space $E$.  Then the cohomology
group $H^i (E, E_0; \Z)$ is zero for $i < l$, and $H^l (E, E_0; \Z)$
contains one and only one cohomology class $u$ whose restriction
\[
	u \mid (F, F_0) \in H^l (F, F_0; \Z)
\]
is equal to the preferred generator $u_F$ for every fiber $F$ of $\xi$.
Furthermore, the correspondence $y \mapsto y \cup u$ maps $H^k (E; \Z)$
isomorphically onto $H^{k + l}(E, E_0; \Z)$.

\end{theorem}

We will be using the Thom isomorphism below in de Rham cohomology, so it will serve
us to re-state the theorem in this context.  The construction below is based on the
work of Schwerdtfeger \cite{schwerdtfeger}.

Let $M$ be a manifold with $E$ an $l$-dimensional vector
bundle over $M$.  Assuming a Euclidean metric on the bundle $E$, let $B$ denote
the unit disk bundle of $E$ (i.e. the set of vectors $v$ with $\| v \| \leq 1$) and
$\partial B = S$ the unit sphere bundle.  We let $\rho$ denote the projection
$\rho : E \rightarrow M$, as
well as its restriction to $B$, $S$, etc.  Let $B_0 = E_0 \cap B$, and let
$\phi : B_0 \rightarrow S$ be the map $v \mapsto \frac{1}{\| v \|}v$.

The cohomology of the pair $(E, E_0)$, then, is studied via forms on $B$
relative to $S$
(i.e. forms on $B$ which vanish on the boundary).  We state the following theorem
from \cite{schwerdtfeger}, noting that we have changed his notation and sign
convention in a consistent and suggestive manner.

% XXXXXXXXXXXXXXXXXXXXXXXXXXXX	THEOREM

\begin{theorem} (Schwerdtfeger)
\label{sch}
Let $E(\Omega) \in \Omega^l(M)$, $\Psi \in \Omega^{l-1}(S)$ be forms which satisfy
\[
	dE(\Omega) = 0, \;\;
	d\Psi = - \rho^\ast(E(\Omega)),	\;\;
	\int\limits_{S_p} \Psi = 1
\]
where $S_p$ denotes the fiber of $S$ over an arbitrary point $p \in M$.  Let
$h : B \rightarrow \R$ be smooth with
\[
	h = \left\{ \begin{array}{ll}
		0	&	\mbox{inside of a $\delta$-neighborhood of the
								zero section}		\\\\
		1	&	\mbox{outside of an $\epsilon$-neighborhood of the
					zero section}	
	\end{array} \right.
\]
where $0 < \delta < \epsilon < 1$.

Then the form
\[
	\tau = \rho^\ast(E(\Omega)) + d(h \cdot \rho^\ast (\Psi)) \in \Omega^l(B, S)
\]
is a representative of the Thom class.

\end{theorem}

Quite clearly, in the case of a manifold,
the forms $E(\Omega)$ and $\Psi$, as defined above, satisfy these conditions,
and hence we may take this to be a definition of the Thom class for $E$ (here again,
we take the restriction of $\Psi$ to the sphere bundle of $E$).  Note that, restricted
to $S$,
\[
	\tau = (\rho^\ast E(\Omega))_{|S} + d(\Psi_{|S}) = 0,
\]
so that $\tau$ vanishes on the boundary.  Note further that on $M$, $\tau = E(\Omega)$.

With this, we have that the map
\[
\begin{array}{rccl}
	\psi	: &	\Omega^k (M)	& \rightarrow	& \Omega^{k+l} (B,S)	
									\\\\
		: &	\omega		& \mapsto	&
					\rho^\ast(\omega) \wedge \tau
\end{array}
\]
induces an isomorphism
\[
	H^k(M) \cong H^{k + l} (E, E_0),
\]
where we are using the natural isomorphisms
$H^{k + l} (E, E_0) \cong H^{k + l}(B, S)$ and $H^k(E) \cong H^k(M)$
induced in both cases by the inclusion of the latter into the former.

% XXXXXXXXXXXXXXXXXXXXXXXXXXXX	subsection: good orbifolds

\subsection{The Case of a Global Quotient}

Let $Q := M/G$ be an $n$-dimensional, oriented, closed global quotient orbifold so that
$M$ is oriented, compact, and smooth and $G$ is finite, and let $E$ be a $G$-bundle of
rank $l$ on $M$ equipped with a Euclidean metric.  Then $E/G$ naturally has the
structure of an orbifold vector bundle
over $Q$.  Say that $E/G$ carries a Euclidean metric (which is precisely a
$G$-invariant metric on $E$).  Recall that a differential form $\omega$ on $Q$
is precisely a $G$-invariant differential form on $M$.  Hence, if $(\Omega^k(M))^G$
denotes the $G$-invariant $k$-forms on $M$, then the quotient projection
\[
	\pi : M \rightarrow M/G = Q
\]
induces
\[
	\pi^\ast : \Omega^k(Q) \rightarrow (\Omega^k (M))^G,
\]
which is clearly an isomorphism of linear spaces.
If we extend $\pi$ to $E \rightarrow E/G$,
then we have a similiar identification
\[
	\Omega^k(E/G) \cong (\Omega^k(E))^G.
\]

Throughout this section, let $\pi$ denote the projection of $M$ onto $Q$ (and its
various extensions to bundles on $M$ and their corresponding orbifold bundles on these $Q$).
With $E$ and $E/G$ the bundles over $M$ and $Q$, respectively (with respective
projections $\rho_M$ and $\rho_Q$), we let $E_0$ and $(E/G)_0$ denote the
collection of nonzero vectors in each of these spaces, $B_M$ and $B_Q$ the ball bundles,
$S_M$ and $S_Q$ the sphere bundles, etc.

Let $u_M \in H^n(E, E_0; \R)$ denote the Thom class of the bundle $E$ over $M$
(tensoring the cohomology group with $\R$), and let $\tau \in \Omega^l(B_M, S_M)$
be the differential form given by Theorem \ref{sch} which represents the Thom class in de Rham
cohomology.
In general, if $\omega \in (\Omega^k(M))^G$ is a $G$-invariant differential form
on $M$, we will identify it with the associated differential form on $Q$.  We denote by
$[\omega]_M$ its class in $H^k(M; \R)$ and $[\omega]_Q$
its class in $H^k(Q ; \R)$ (and respectively, on $E$, $(E, E_0)$, etc.).  So in
this notation, $u_M = [\tau]_M$.

We have that the map
\[
\begin{array}{rccl}
	\psi : &	H^k(M; \R)	& \rightarrow	& H^{k+l}(E,E_0;\R )	
									\\\\
		 : &	[\omega]_M	& \mapsto	&
					[\rho_M^\ast(\omega) \wedge \tau_M]_M
\end{array}
\]
is an isomorphism (again useing the canonical isomorphisms
$H^k(M;\R) \cong H^k(E;\R)$ and $H^{k+l}(E, E_0; \R) \cong H^{k+l}(B_M, S_M;\R)$).

Note that $\tau$ is $G$-invariant, as it is defined in terms of the forms $\Psi$
and $E(\Omega)$, which are $G$-invariant whenever the metric is, and the function $h$
can clearly be chosen to be $G$-invariant.  Hence,
$\tau \in \Omega^l(B_Q, S_Q)$ (via its identification with
$(\Omega^l(B,S))^G$), and as $d\tau = 0$ clearly, $\tau$ represents a cohomology
class in $H^l(Q; \R)$.

Consider the map $\hat{\psi} : H^k(Q; \R) \rightarrow H^{k+l}(E/G, (E/G)_0; \R)$,
where $\hat{\psi}([\omega]_Q) = [\rho_Q^\ast(\omega) \wedge \tau]_Q$.  We will
represent this map using the isomorphisms $\Omega^k(E/G) \cong (\Omega^k(E))^G$ and
$\Omega^k(B_Q, S_Q) \cong (\Omega^k(B_M, S_M))^G$ given above, and the obvious
isomorphism
\[
	H^{k+l}(E/G, (E/G)_0; \R) \cong H^{k+l}(B_Q, S_Q; \R).
\]
Hence, this map can be expressed on forms as
$\hat{\psi}(\omega) = \rho_M^\ast(\omega) \wedge \tau$
for $\omega \in (\Omega^k(M))^G$.

Before we deal with $\hat{\psi}$, however, we need a lemma which will help us
relate $\Omega^k(M)$ and $(\Omega^k(M))^G$.  Essentially, it claims that if a
$G$-invariant form is exact in $\Omega^k(M)$, then it is exact in $(\Omega^k(M))^G$.

% XXXXXXXXXXXXXXXXXXXXXXXXXXXX	LEMMA

\begin{lemma}
\label{gexact}

Suppose that $\eta_1 \in \Omega^{k-1}(M)$ with $d\eta_1 \in (\Omega^k(M))^G$.
Then there is an $\eta_2 \in (\Omega^{k-1}(M))^G$ with $d\eta_1 = d\eta_2$.

\end{lemma}

\noindent	{\it Proof:}

Using the averaging map, set
\[
	\eta_2 := \frac{1}{|G|} \sum\limits_{\gamma \in G} \gamma^\ast \eta_1 .
\]
Then 
\[
\begin{array}{rcl}
	d\eta_2		&=&	d \left( \frac{1}{|G|} \sum\limits_{\gamma \in G} \gamma^\ast \eta_1 \right)		\\\\
			&=&	\frac{1}{|G|} \sum\limits_{\gamma \in G} \gamma^\ast d\eta_1				\\\\
			&=&	\frac{1}{|G|} \sum\limits_{\gamma \in G} d\eta_1					\\\\
			&&	\mbox{(as $d\eta_1$ is $G$-invariant)}						\\\\
			&=&	d\eta_1 .
\end{array}
\]
Moreover, for each $\gamma_0 \in G$,
\[
\begin{array}{rcl}
	\gamma_0^\ast \eta_2	&=&	\gamma_0^\ast \left( \frac{1}{|G|} \sum\limits_{\gamma \in G} \gamma^\ast d\eta_1 \right)	\\\\
			&=&	\frac{1}{|G|} \sum\limits_{\gamma \in G} \gamma_0^\ast \gamma^\ast d\eta_1			\\\\
			&=&	\frac{1}{|G|} \sum\limits_{\gamma \in G} (\gamma\gamma_0)^\ast d\eta_1				\\\\
			&=&	\frac{1}{|G|} \sum\limits_{\gamma \in G} \gamma^\ast d\eta_1				\\\\
			&&	\mbox{(as multiplication by an element simply permutes}				\\
			&&	\mbox{the elements of a group)}							\\\\
			&=&	\eta_2,
\end{array}
\]
so that $\eta_2$ is $G$-invariant.

\bigskip \noindent Q.E.D. \bigskip

% XXXXXXXXXXXXXXXXXXXXXXXXXXXX	CLAIM

\begin{claim}

The map $\hat{\psi} : H^k(Q; \R) \rightarrow H^{k+l}(E/G, (E/G)_0; \R)$
defined above is injective.

\end{claim}

\noindent	{\it Proof:}

Suppose that $\omega_1, \omega_2 \in (\Omega^k(M))^G$ represent classes in
$H^k(E/G; \R)$ such that $\hat{\psi}([\omega_1]_Q) = \hat{\psi}([\omega_2]_Q)$.
Hence, $[\rho_Q^\ast(\omega_1) \wedge \tau]_Q = [\rho_Q^\ast(\omega_2) \wedge \tau]_Q$.

As $\omega_1$ and $\omega_2$ are closed, they represent classes
$[\omega_1]_M, [\omega_2]_M \in H^k(M; \R)$. So as $\psi$ is known to be
an isomorphism here, and as $\psi([\omega_1]_M) = \psi([\omega_2]_M)$,
we have that $[\omega_1]_M = [\omega_2]_M$.  So there is an
$\eta \in \Omega^{k-1}(M)$
such that
\[
	\omega_1 = \omega_2 + d\eta .
\]
Note that, as
\[
	d\eta = \omega_1 - \omega_2 ,
\]
and the $\omega_i$ are $G$-invariant, $d\eta$ is $G$-invariant.  Applying Lemma \ref{gexact},
we can take $\eta$ to be $G$-invariant, and hence $[\omega_1]_Q = [\omega_2]_Q$.  So $\hat{\psi}$ is injective.

\bigskip \noindent Q.E.D. \bigskip

% XXXXXXXXXXXXXXXXXXXXXXXXXXXX	CLAIM
\begin{claim}

The map $\hat{\psi} : H^k(Q; \R) \rightarrow H^{k+l}(E/G, (E/G)_0; \R)$
is surjective.

\end{claim}

\noindent	{\it Proof:}

Let $\zeta \in (\Omega^{k+l}(E, E_0))^G$ be a closed form.
Then $\zeta$ represents a class $[\zeta]_M \in H^{k+l}(E, E_0; \R)$,
so that as $\psi$ is an isomorphism, there is a closed $\omega \in \Omega^k(M)$ and an
$\eta \in \Omega^{k-1}(M)$ such that
$\zeta = \rho_M^\ast(\omega) \wedge \tau + d\eta$.  As $\omega$ is closed,
$d\omega = 0$ is $G$-invariant, so that $\omega$ can be taken to be
$G$-invariant.  So $\omega$ represents a class in $H^k(Q; \R)$.
Moreover, $d\eta = \zeta - \rho_M^\ast(\omega) \wedge \tau$
is $G$-invariant, so that $\eta$ can be taken to be,
Therefore, $\zeta$ is cohomologous to $\rho_M^\ast(\omega) \wedge \tau$
in $H^{n+l}(Q; \R)$.

\bigskip \noindent Q.E.D. \bigskip

% XXXXXXXXXXXXXXXXXXXXXXXXXXXX	CLAIM
\begin{claim}

The cohomology class $[\tau]_Q \in H^l(E/G, (E/G)_0)$ does not depend on
the connection.

\end{claim}

\noindent	{\it Proof:}

We need only show that the map $\pi^\ast : H^l (E/G, (E/G)_0) \rightarrow H^l(E, E_0)$
induced by the quotient projection $\pi : E \rightarrow E/G$ is injective.
for the cohomology class of $\pi^\ast \tau$ is unique by Theorem \ref{milnorthom}.
So suppose $\pi^\ast \omega_1 - \pi^\ast \omega_2 = d\eta$ for closed
$\omega_1, \omega_2 \in \Omega^k (B_Q, S_Q)$ and some $\eta \in \Omega^{k-1}(B_M, S_M)$.
Then as $d\eta = \pi^\ast (\omega_1 - \omega_2)$ is $G$-invariant, $\eta$ can be chosen
to be $G$-invariant, so that $\eta$ is in the image of $\pi^\ast$.  Hence $\omega_1 - \omega_2$
is exact.

\bigskip \noindent Q.E.D. \bigskip

Hence, we see that $\tau$ represents a unique cohomology class $u_Q \in H^l(E/G, (E/G)_0; \R)$
which induces an isomorphism
$\rho_Q^\ast (\cdot) \cup u_Q : H^k(Q; \R) \rightarrow H^{k+l}(E/G,(E/G)_0;\R )$.

In summary, we state:

\begin{proposition}[The Thom Isomorphism for Global Quotients]
Let $Q = M/G$ be an $n$-dimensional, oriented, closed global quotient orbifold so that
$M$ is oriented, compact, and smooth and $G$ is finite, and let $E$ be a $G$-bundle
on $M$.  Then, with the form $\tau$ as defined above, the map
\[
\begin{array}{rccl}
	\hat{\psi}
		: &	H^k(Q; \R)	& \rightarrow	& H^{k+l}(E/G,(E/G)_0;\R )
									\\\\
		: &	[\omega]_Q	& \mapsto	&
					[\rho_Q^\ast(\omega) \wedge \tau ]_Q
\end{array}
\]
is an isomorphism.  Hence, if $u_Q = [\tau]_Q$ denotes the cohomology class
in $H^l (E/G, (E/G)_0 ; \R)$ represented by $\tau$, then the map
$\rho_Q^\ast(\cdot) \cup u_Q : H^k(Q; \R) \rightarrow H^{k+n}(E/G,(E/G)_0;\R )$ is an isomorphism.
Moreover, $u_Q$ does not depend on the metric or connection.
\end{proposition}

In the next section, we will apply these results locally to a bundle $E$
over a general orbifold by considering each set $U$ with chart
$\{ V, G, \pi \}$ to be a global quotient $V/G$, and the bundle $E$ to
be the quotient $(\pi^\ast E)/G$.
However, in these cases, the set $V$ is not
compact.  In order to justify such an application, we note that
although Schwerdtfeger's results only apply to compact manifolds, Theorem
\ref{milnorthom} applies in general, so that there is a Thom class in
$H^k (\pi^\ast (E), (\pi^\ast (E))_0)$ with the desired properties.
Moreover, as we can always take $V$ to be a bounded open ball
in $\R^n$ and extend the metric and connection to a larger ball in $\R^n$
which contains $\overline{V}$.  Hence, as the Thom class can be characterized
by a completely local property (that its restriction to each fiber $F$ is the
preferred generator of $H^k(F, F_0)$ with respect to the orientation), the
class of $\tau$ must coincide with the Thom class on the interior of this
ball.  So in this case, the above construction is still valid.

% XXXXXXXXXXXXXXXXXXXXXXXXXXXX	subsection: general orbifold

\subsection{The Case of a General Orbifold}

For the case of a general closed orbifold $Q$ with orbibundle
$\rho: E \rightarrow Q$, we will
use the fact that $Q$ is locally a global quotient; i.e. that each point is
contained in a neighborhood which is a global quotient.  The proof of the Thom
Isomorphism Theorem, then, will involve
an induction following Milnor and Stasheff \cite{milnor}.  We note that for each open
$U \subset Q$ that is uniformized by $\{ V, G, \pi \}$, $U$
is given the structure of a global quotient orbifold via $\tilde{\pi} : V / G \cong U$.
Hence, applying the note in the last subsection, the Thom Isomorphism Theorem is known locally.

First, we let $E$ carry a Euclidean metric.  Then we may take the definition
of the Thom class, $\tau$, as given above, in terms of the global forms $\Psi$
and $E(\Omega)$.  As $d\tau = 0$, $\tau$ represents a cohomology class
$[\tau] =: u \in H^n(B, S; \R) \cong H^n(E, E_0; \R)$, which we define to
be the Thom class of $Q$ (here, $B$ and $S$ are defined to be the ball
bundle and sphere bundle, respectively, of $E$ with respect to the metric,
as in the last subsection).  Note that, for each uniformized $U \subset Q$,
we may give $U$ the restricted metric of $Q$, and then $\tau_{\mid U}$ is
the Thom form of $U$.  We now proceed with the induction.

Suppose $Q = U_1 \cup U_2$ is the union of two open sets, with $U_i$ having
chart $\{ V_i, G_i, \pi_i \}$ for $i = 1, 2$, such that each $V_i$ is a bounded
ball in $\R^n$, and each $G_i$ acts linearly.  Then applying the note above and the case
for global quotients, the Thom Isomorphism holds for each $E_{\mid U_i}$.
Note further that, as $W := U_1 \cap U_2 \subset U_1$, $W$ is also given the structure
of a quotient via the uniformization of $U_1$.

We have the following two Meier-Vietoris sequences (using coefficients in $\R$ throughout):
\[
	\cdots
	\rightarrow
	H^{k - 1}(W)
	\rightarrow
	H^k(Q)
	\rightarrow
	H^k(U_1)
	\oplus H^k(U_2)
	\rightarrow
	H^k(W)
	\rightarrow
	\cdots
\]
and
\[
	\cdots
	\rightarrow
	H^{k + l - 1}(E_{|W}, (E_{|W})_0)
	\rightarrow
	H^{k + l}(E, E_0)
	\rightarrow
\]\[
	\rightarrow
	H^{k + l}(E_{\mid U_1}, (E_{\mid U_1})_0)
	\oplus H^{k + l}(E_{\mid U_2}, (E_{\mid U_2})_0)
	\rightarrow
	H^{k + l}(E_{|W}, (E_{|W})_0)
	\rightarrow
	\cdots
\]
Using the fact that the isomorphism is known for all but one step in this sequence,
we have the following for each $k$,
\[
\begin{array}{cccc}
	H^{k - 1}(W)	&
	\rightarrow	&
	H^k(Q)		&
	\rightarrow		\\\\
	\downarrow\cong	&&
	\downarrow\ast	&	\\\\
	H^{k + l - 1}(E_{|W}, (E_{|W})_0)
			&
	\rightarrow	&
	H^{k + l}(E, E_0)
			&
	\rightarrow
\end{array}
\]\[
\begin{array}{cccc}
	\rightarrow	&
	H^k(U_1)\oplus H^k(U_2)
			&
	\rightarrow	&
	H^k(W)			\\\\
			&
	\downarrow\cong	&
			&
	\downarrow\cong		\\\\	
	\rightarrow	&
	H^{k + l}(E_{\mid U_1}, (E_{\mid U_1})_0)
	\oplus H^{k + l}(E_{\mid U_2}, (E_{\mid U_2})_0)
			&
	\rightarrow	&
	H^{k + l}(E_{|W}, (E_{|W})_0)
\end{array}
\]
Each of the vertical isomorphisms is given by $\rho^\ast (\cdot) \wedge \tau$ on the
level of forms; similarly, by $\rho^\ast (\cdot) \cup u$ in cohomology (with $\tau$ and $u$
appropriately restricted), so that the diagram clearly commutes.  Applying the
Five Lemma gives us that $(\ast)$ is also an isomorphism, and
that the isomorphism is given by $\rho^\ast (\cdot) \wedge \tau$.  Hence, we have shown
the Thom isomorphism in this case.

Now, suppose $Q$ is any closed orbifold with orbifold vector bundle $E$ of rank
$l$, and let $\{ U_i \}_{i=1}^k$ be a cover
of $Q$ such that each $U_i$ is uniformized by
$\{ V_i, G_i, \pi_i  \}$.  If $k=1$, then $Q$ is a global quotient,
and Thom isomorphism
is known.  For $k > 1$, assuming as our inductive hypothesis that the Thom
isomorphism holds for
$E_{\mid U_1 \cup \cdots \cup U_{k-1}}$, applying the above argument to
the two sets $U_1 \cup \cdots \cup U_{k-1}$ and $U_k$ shows that it holds for $E$.

With this, we have proven the following.
\begin{theorem}[The Thom Isomorphism for a Closed Orbifold]
\label{thomgeneral}
Let $Q$ be a closed oriented orbifold and $\rho:E \rightarrow Q$
an oriented orbifold vector bundle of rank $l$ over $Q$.  Let $\tau$
be defined as in Theorem \ref{sch}, and then the map
\[
\begin{array}{rccl}
	\hat{\psi}
		: &	H^k(Q; \R)	& \rightarrow	& H^{k+l}(E,E_0;\R )
									\\\\
		: &	[\omega]_Q	& \mapsto	&
					[\rho_Q^\ast(\omega) \wedge \tau ]_Q
\end{array}
\]
is an isomorphism.  Moreover, with respect to each chart $\{ V, G, \pi \}$,
$\pi^\ast \tau$ restricts to the preferred
generator of $H^l(F, F_0 \; \Z)$ (using the isomorphism between
$H^l(F, F_0 \; \R)$ and $H^l(F, F_0 \; \Z) \otimes \R$).
In particular, the class of $\tau$ does not depend on the metric or connection.
\end{theorem}

Note that, by its construction, it is trivial that the restriction of
$\tau$ to $Q$ is $E(\Omega)$ (and hence the restriction of $u$ to $Q$
is the cohomology class of $E(\Omega)$ in $H^l (Q)$).  Therefore,
$E(\Omega)$ is closed, and cohomology class it represents in $H^l(Q)$
does not depend on the metric.  We refer to this class as $e(E)$,
the {\bf Euler class of $E$ as an orbifold bundle}.  In the case
$E = TQ$, $e(TQ) =: e(Q)$ is simply the {\bf Euler class of $Q$ as
an orbifold.}  In general it does not coincide with the Euler class
of the underlying space of $Q$, but is a rational multiple thereof.

% XXXXXXXXXXXXXXXXXXXXXXXXXXXXXXXXXXXXXXXXXXXXXXXXXXXXXXXXXXXXXXXXXXX
%			Section: the theorems on the level of cohomology
% XXXXXXXXXXXXXXXXXXXXXXXXXXXXXXXXXXXXXXXXXXXXXXXXXXXXXXXXXXXXXXXXXXX

\section{The Theorems on the Level of Cohomology}
\label{thrmcohom}

\subsection{Invariance of the Integrands on the Metric}

We return to the case where $Q$ has boundary,
$M = \partial Q$, and $E = TQ$.  We have the formula
(Proposition \ref{phdependent})
\[
	\mbox{ind}(X)		=
	\left\{	\begin{array}{ll}
		\chi_{orb}^\prime (Q) + \int\limits_M \alpha^\ast(\Psi),
			&			n = 2m,				\\\\
		\int\limits_M \alpha^\ast(\Psi),
			&			n = 2m + 1,
	\end{array} \right.
\]
which generalizes the Poincar\'{e}-Hopf Theorem to the case of a compact
orbifold with boundary.  In this section, we characterize the cohomology class
of $\Psi$ in $H^n (STQ_M)$ in order to show that it does not depend on the metric.
This section again follows \cite{sha}.

Recall that $d\Psi = - \rho^\ast E(\Omega)$ on $STQ$.  However, as
$TQ_{|M}$ is isomorphic to
$TM \oplus \nu$ where $\nu$ denotes the trivial bundle on $M$ of rank 1,
and hence that $\Omega_{i, j} = 0$ whenever $i > 1$ or $j > 1$, $E(\Omega)$
clearly vanishes on $M$.  So $\Psi$ is a closed
form on $STQ_{|M}$, and hence represents a cohomology class
$\Upsilon \in H^n (STQ_{|M}; \R)$.  We have seen that
\[
	\alpha_0^{\ast} (\Upsilon ) = \frac{1}{2}e(TM),
\]
where $\alpha_0$ is the outward-pointing unit normal vector field on $M$
and $e(TM)$ denotes the cohomology class of $E(\Omega_{| M})$ in $H^{n-1}(M; \R)$.
Moreover, the integral of $\Psi$ on a fiber of $STQ_{|M}$ over a point $p$
is $\frac{1}{|I_p|}$.  Thus, the following is clear.

% XXXXXXXXXXXXXXXXXXXXXXXXXXXX	CLAIM

\begin{claim}

Let $p \in U_p \subset M$ where $U_p$ is uniformized by $\{ V_p, G_p, \pi_p \}$.
If $\iota : S^n \rightarrow \pi_p^{\ast} STQ_{|M}$ is an orientation-preserving
isometry from $S^n$ to a fiber of the sphere bundle of $\pi_p^{\ast}TQ$
(i.e. the sphere bundle $STQ$ pulled back over $V$), then
then $\iota^{\ast} \pi_p^{\ast} (\Upsilon) = \frac{1}{|G_p|}s^n$, where $s^n$ denotes
the canonical generator of $H^n(S_n; \Z)$.

\end{claim}

Note that in the case where $G = 1$, we can identify $U$ with $V$ via $\pi_p$, so that
$\iota^{\ast} \pi_p^{\ast} \Upsilon = \iota^{\ast} \Upsilon$.

With this, we may characterize the class $\Upsilon$ by developing the Gysin sequence of the
tangent bundle.  Let $E$ now denote the restriction
$TQ_{|M}$ and $E_0$ the nonzero vectors of $E$.  We follow Milnor and Stasheff \cite{milnor}
and Sha \cite{sha}.  All coefficient groups are understood to be $\R$.

To simplify the notation, we let $E$ denote $TQ_{|M}$.
We have the cohomology exact sequence
\[
	\cdots \longrightarrow H^j (E, E_0) \longrightarrow H^j(E)
	\longrightarrow H^j(E_0) \longrightarrow H^{j+1} (E, E_0)
	\longrightarrow \cdots
\]
We may replace $H^j(E)$ with $H^j(M)$ using the natural isomorphism induced by $\rho$.
Similarly, applying the Thom Isomorphism, we may replace $H^j(E, E_0)$ with
$H^{j - n - 1}(M)$.  The restriction map
$H^j (E, E_0) \rightarrow H^j(E)$ composed with $\cdot \cup u$ results in
$\cdot \cup e$ (where $e$ is the Euler class of the orbifold bundle $E$):
\[
	\cdots \longrightarrow H^{j - n - 1} (M)
	\stackrel{\cup e}{\longrightarrow} H^j(M)
	\longrightarrow H^j(E_0) \longrightarrow H^{j - n} (M)
	\stackrel{\cup e}{\longrightarrow} \cdots
\]
Now $H^j(E_0)$ is canonically isomorphic to $H^j(SE)$, so that again using $\rho$
for the restricted projection $\rho : SE \rightarrow M$ we obtain
\[
	\cdots \longrightarrow H^{j - n - 1} (M)
	\stackrel{\cup e}{\longrightarrow} H^j(M)
	\stackrel{\rho^{\ast}}{\longrightarrow} H^j(SE)
	\longrightarrow H^{j - n} (M)
	\stackrel{\cup e}{\longrightarrow} \cdots
\]
However, we have noted that the Euler class $e$ of $E$ is zero, so that setting
$j = n$, we obtain
\[
	0 \longrightarrow H^{n} (M)
	\longrightarrow H^n(M)
	\stackrel{\rho^{\ast}}{\longrightarrow} H^n(STQ)
	\longrightarrow H^{0} (M)
	\stackrel{\cup e}{\longrightarrow} 0
\]
Finally, we choose $\iota :S^n \rightarrow SE_p$ to be an isometry (as above)
to a fiber over a point $p$ with trivial isotropy, so that 
\[
\begin{array}{ccccccccc}
	0	&	\longrightarrow
		&	H^n(M)
		&	\stackrel{\rho^\ast}{\longrightarrow}
		&	H^n(SE_{\mid M})
		&	\rightarrow
		&	H^0(M)
		&	\stackrel{\cup e}{\longrightarrow}
		&	0					\\\\
		&&&&	\downarrow	\iota^\ast
		&&	\downarrow	\cong
		&&						\\\\
		&&&&	H^n(S^n)
		&	\stackrel{\cong}{\longrightarrow}
		&	H^0(\mbox{point})
		&&
\end{array}
\]
gives the split exact sequence
\[
	0 \longrightarrow H^{n} (M)
	\stackrel{\rho^{\ast}}{\longrightarrow} H^n(SE)
	\stackrel{\iota^{\ast}}{\longrightarrow} H^n(S^n)
	\longrightarrow	0
\]
With this, as $\alpha_0^\ast$ is a left inverse of $\rho^\ast$
(recall that $\alpha_0: M \rightarrow SE$ is the section induced by
the outward unit normal vector field on $M$), we have
$H^n(SE) = \rho^\ast(H^n(M)) \oplus (\alpha_0^\ast)^{-1}(0)$.
With respect to this decomposition, based on the properties of
$\Psi$, $\Upsilon$ factors into
$\frac{1}{2} \rho^\ast(e) + \iota_{\mid (\alpha_0^\ast)^{-1}(0)}^{-1}(s^n)$.
Hence, this characterizes $\Upsilon$, and in particular shows that it
does not depend on the choices made in the definition of $\Psi$.
Note that choosing $\iota$ to map to the fiber over a singular point $p$
will introduce a coefficient of $\frac{1}{|G_p|}$ in the above expression.

% XXXXXXXXXXXXXXXXXXXXXXXXXXXX	PROPOSITION

With this, Proposition \ref{phdependent} becomes
\begin{theorem}[The First Poincar\'{e}-Hopf Theorem for Orbifolds with Boundary]
\label{pcb1}
Let $Q$ be a compact orbifold with  boundary $M$.  Let
$X$ be vector field on $Q$ which has a finite number of singularities,
all of which occuring on the interior of $Q$.  Then
\[
	\mbox{ind}(X)		=
	\left\{	\begin{array}{ll}
		\chi_{orb}^\prime (Q) +	\alpha^\ast \Upsilon ([M]),
			&			n = 2m,				\\\\
		\alpha^\ast \Upsilon ([M]),
			&			n = 2m + 1.
	\end{array} \right.
\]
\end{theorem}

%% file: thesisch4.tex
% Thesis Chapter 4, Christopher Seaton, seatonc@colorado.edu
% XXXXXXXXXXXXXXXXXXXXXXXXXXXXXXXXXXXXXXXXXXXXXXXXXXXXXXXXXXXXXXXXXXX
% XXXXXXXXXXXXXXXXXXXXXXXXXXXXXXXXXXXXXXXXXXXXXXXXXXXXXXXXXXXXXXXXXXX
%			Chapter 4: Gauss-Bonnet in C-R Cohomology
% XXXXXXXXXXXXXXXXXXXXXXXXXXXXXXXXXXXXXXXXXXXXXXXXXXXXXXXXXXXXXXXXXXX
% XXXXXXXXXXXXXXXXXXXXXXXXXXXXXXXXXXXXXXXXXXXXXXXXXXXXXXXXXXXXXXXXXXX
%
%	section 1:	Introduction for this chapter
%	section 2:	Summary of chen-ruan orbifold cohomology
%	section 3:	C-R cohomology for orbifolds with boundary
%	section 4:	The GB and PH theorems for closed orbifolds
%	section 5:	The GB and PH theorems for orbifolds with boundary
%

\chapter{The Gauss-Bonnet Integrand in Chen-Ruan Cohomology}
\label{chcr}

% XXXXXXXXXXXXXXXXXXXXXXXXXXXXXXXXXXXXXXXXXXXXXXXXXXXXXXXXXXXXXXXXXXX
%			Section: Intro this chapter
% XXXXXXXXXXXXXXXXXXXXXXXXXXXXXXXXXXXXXXXXXXXXXXXXXXXXXXXXXXXXXXXXXXX

\section{Introduction}

The goal of this chapter is to examine the results of
the previous chapter in terms of the orbifold cohomology developed in
\cite{chenruan}.  In particular, we are interested in cohomology classes
corresponding to those of $E(\Omega)$ and $\Psi$.  Roughly speaking,
the Chen-Ruan orbifold cohomology of an orbifold $Q$ contains the usual cohomology
of $Q$ as a direct summand, but contains as well the cohomology groups
of the twisted sectors, corresponding to irreducible components of the
singular set of $Q$.  With respect to this decomposition, the new
characteristic classes will project to the usual ones.  They will,
however, have additional lower-degree terms, which correspond to the
contributions of the singular sets.  Our results, then, will involve
the Euler characteristic of the underlying topological space
instead of the orbifold Euler characteristic.
Taking the point of view that an orbifold structure is a generalization
of a differentiable structure on a manifold, we obtain results for
orbifolds much more in keeping with the original Gauss-Bonnet
and Poincar\'{e}-Hopf Theorems.

Throughout this chapter, orbifolds and orbifolds with boundary
will be taken to admit almost complex structures.

% XXXXXXXXXXXXXXXXXXXXXXXXXXXXXXXXXXXXXXXXXXXXXXXXXXXXXXXXXXXXXXXXXXX
%			Section: Review Ch-R Cohomology
% XXXXXXXXXXXXXXXXXXXXXXXXXXXXXXXXXXXXXXXXXXXXXXXXXXXXXXXXXXXXXXXXXXX

\section{Chen-Ruan Orbifold Cohomology}

In this chapter, we will be working with the orbifold cohomology theory
developed by \cite{chenruan}.
We will not develop this cohomology theory here, but will collect a summary for the sake of making the
notation explicit.  For the most part, we follow the notation in 
\cite{chenruan}, \cite{ruangwt}, and \cite{ruansgt}.

Let $Q$ be an orbifold, and select for each $p \in Q$ a chart $\{ V_p, G_p, \pi_p \}$ at $p$.
Then the set
\[
	\tilde{Q} = \{ (p, (g)_{G_p} ) : p \in Q, g \in G_p \}
\]
(where $(g)_{G_p}$ is the conjugacy class of $g$ in $G_p$) is naturally an orbifold, with local charts
\[
	\{ \pi_{p, g} : (V_p^g, C(g) ) \rightarrow V_p^g / C(g) : p \in Q, g \in G_p \} ,
\]
where $V_p^g$ is the fixed point set of $g$ in $V_p$ and $C(g)$ is the centralizer of $g$ in $G_p$.
If $Q$ is closed, then $\tilde{Q}$ is closed, but it need not be connected, and its connected
components need not be of the same dimension.  An equivalence relation can be placed on the
elements of the groups $G_p$ so that if $T$ denotes the set of equivalence classes and $(g)$ the
equivalence class of $g$,
\[
	\tilde{Q} = \bigsqcup\limits_{(g) \in T } \tilde{Q}_{(g)}
\]
where
\[
	\tilde{Q}_{(g)} = \{ (p, (g^\prime)_{G_p}) : g^\prime \in G_p , (g^\prime)_{G_p} \in (g) \} .
\]
The map $\pi : \tilde{Q} \rightarrow Q$ with $(p, (g)_{G_p}) \mapsto p$, is a $C^\infty$
map.

If $Q$ is an almost complex orbifold, a function $\iota : \tilde{Q} \rightarrow \Q$ is
defined which is constant on the connected components of $\tilde{Q}$.
If $n_{(g)}$ denotes the codimension of $\tilde{Q}_{(g)}$ in $Q$, then
$2\iota_{(g)} \leq n_g$, with equality only when $g = 1$.  This is
called the degree shifting number of $(g)$.
The orbifold cohomology groups are defined by
\[
	H_{orb}^d (Q)	=	\bigoplus_{(g) \in T} H^{d - 2\iota_{(g)}} (\tilde{Q}_{(g)}),
\]
where the groups on the right side are the usual de Rham cohomology groups of the orbifolds
$\tilde{Q}_{(g)}$.

Since each $\tilde{Q}_{(g)}$ can be realized as a subset of $Q$, geometric constructions
(i.e. bundles and their sections) on $Q$ can be naturally extended to geometric
constructions on $\tilde{Q}$.  In what follows, we wish to extend the characteristic
classes of bundles over $Q$ to characteristic classes of associated bundles over $\tilde{Q}$;
however, pulling back such bundles via $\pi$ will be insufficient.  In particular, if
$E \rightarrow Q$ is a rank $k$ orbibundle and $p \in Q$ is a singular point contained in a singular
set of dimension $l < k$, then the maximal vector space in a fiber over $p$ has dimension $l$.
This implies that any $k$-form on $Q$ is zero on $p$ (recall that it is required of sections $s$ of
orbibundles that for each $q \in Q$, $s(q)$ is contained in the subspace $E^q$ of the fiber over $q$ which is
fixed by $G_q$).  In particular, any form
representing the Euler class of $E$ is zero at $p$, so that the pull-back $\pi^{\ast}$ of
this form will be zero on each connected component $\tilde{Q}_{(g)}$ of $\tilde{Q}$
of dimension is less than $k$.  This is particularly disappointing in the case of the tangent
bundle, in that no twisted sectors make contributions to the Euler class.

Instead, we will associate to each bundle $E \rightarrow Q$ a bundle $\tilde{E} \rightarrow \tilde{Q}$
whose dimension on each component $\tilde{Q}_{(g)}$ is equal to $k - n + l$, where $k$ is the rank of $E$,
$n$ is the dimension of $Q$, and $l$ is the dimension of $\tilde{Q}_{(g)}$ (i.e. the rank of $E$
minus the codimension of $\tilde{Q}_{(g)}$ in $Q$).  We will then apply the Chern-Weil construction
to a connection on $\tilde{E}$ in order to define characteristic classes in $H_{orb}^{\ast} (Q)$ which are
invariants of $E$.

% XXXXXXXXXXXXXXXXXXXXXXXXXXXXXXXXXXXXXXXXXXXXXXXXXXXXXXXXXXXXXXXXXXX
%			Section: Ch-R Cohomology with Boundary
% XXXXXXXXXXXXXXXXXXXXXXXXXXXXXXXXXXXXXXXXXXXXXXXXXXXXXXXXXXXXXXXXXXX

\section{Chen-Ruan Orbifold Cohomology for Orbifolds with Boundary}

In this section, we generalize orbifold cohomology to the case of orbifolds
with boundary.  This is a straightforward generalization following \cite{chenruan}.

Let $Q$ be an $n$-dimensional orbifold with boundary $M$.  Again, we let
\[
	\tilde{Q} = \{ (p, (g)_{G_p} ) : p \in Q, g \in G_p \} .
\]
Then we have:

% XXXXXXXXXXXXXXXXXXXXXXXXXXXX  LEMMA

\begin{lemma}
The set $\tilde{Q}$ is naturally an orbifold with boundary, with projections
given by
\[
	\{ \pi_{p, g} : (V_p^g, C(g) ) \rightarrow V_p^g / C(g) : p \in Q, g \in G_p \} ,
\]
for each chart $(V_p, G_p, \pi_p)$ at $p \in Q$.
Here, $V_p^g$ denotes the fixed point set of $g$ in $V_p$ and $C(g)$ is the
centralizer of $g$ in $G_p$.  If $Q$ is compact, then so is $\tilde{Q}$.  The map
$\pi : \tilde{Q} \rightarrow Q$ is a ${\mathcal C}^\infty$ map.
\end{lemma}

For the proof of this lemma for the case that $Q$ does not have boundary, see
\cite{chenruan}.  The proof for the case with boundaries is identical.

% XXXXXXXXXXXXXXXXXXXXXXXXXXXX  LEMMA

\begin{lemma}

Let $Q$ be an orbifold with boundary $M$.  Then $\partial \tilde{Q} = \tilde{M}$.

\end{lemma}

\noindent       {\it Proof:}

Let $(p, (g)_{G_p})$ be a point in $\partial \tilde{Q}$.  Then $(p, (g)_{G_p})$
is contained in a chart of the form $\{ V_p^g, C(g), \pi_{p, g} \}$, induced
by a chart $\{ V_p, G_p, \pi_p \}$ for $Q$ at $p$.  As $(p, (g)_{G_p})$ is in
th boundary of $\tilde{Q}$, $V_p^g$ is diffeomorphic to $\R_+^k$ for some $k$.
However, as $G_p < O(n)$, its fixed-point set in $\R^n \supset V_p$
is a subspace, so that
$V_p$ must be diffeomorphic to $\R_+^n$ (where $n$ is the dimension of $Q$).
Therefore, $p \in M$.

Conversely, suppose $p \in M$ is a point in the boundary of $Q$, and let
$\{ V_p, G_p, \pi_p \}$ be a chart at $p$.  Then $V_p \cong \R_+^n$, and any
lift $\tilde{p}$ of $p$ into $V_p$ is contained in $\partial V_p$.  For any
$g \in G_p$, the element $(p, (g)_{G_p})$ of $\tilde{Q}$ is covered by the chart
$\{ V_p^g, C(g), \pi_{p, g} \}$, and any lift of $(p, (g)_{G_p})$ into $V_p^g$
is clearly an element of $\partial V_p^g$. Therefore, $(p, (g)_{G_p})$ represents
a point in $\partial \tilde{Q}$.  With this, we note that any point
$(p, (g)_{G_p}) \in \tilde{M}$ arises in such a way, and is contained in a chart
for $\tilde{M}$ induced by a chart for $\tilde{Q}$ (and hence by a chart for $Q$).

\bigskip \noindent Q.E.D. \bigskip

We review the description of the connected components of $\tilde{Q}$, treating the
case that $Q$ has boundary.

Let $\{ V_p, G_p, \pi_p \}$ be an orbifold chart for $Q$ at a point $p \in Q$,
and let $q \in U_p = \pi_p (V_p)$.  Let $\{ V_q, G_q, \pi_q \}$ be an orbifold
chart at $q$ with $U_q \subset U_p$, and then the definition of an orbifold gives
us an injection
$\lambda_{qp} : \{ V_q, G_q, \pi_q \} \rightarrow \{ V_p, G_p, \pi_p \}$.
The injective homomorphism $f_{qp} : G_q \rightarrow G_p$ is well-defined up
to conjugation, so that it defines for each conjugacy class $(g)_{G_q}$
a conjugacy class $(f_{qp}(g))_{G_p}$.  We say that
$(g)_{G_q} \sim (f_{qp}(g))_{G_p}$, which defines an equivalence relation on the
elements of the local groups.  Let $(g)$ denote the
equivalence class of a group element $g$; note that it
is no longer important to state the particular local group from which $g$ was
taken.  For each equivalence class, we let
\[
	\tilde{Q}_{(g)} := \{ (p, (h)_{G_p}) | h \in G_p, (h) =  (g) \} ,
\]
and then $\tilde{Q} = \bigsqcup\limits_{(g) \in T} \tilde{Q}_{(g)}$.  In particular,
following \cite{chenruan}, we call $\tilde{Q}_{(1)}$ the {\bf nontwisted sector}
and each $\tilde{Q}_{(g)}$ for $g \neq 1$ a {\bf twisted sector}.  It is worth
noting that in the case that $Q = M/G$ with $M$ a manifold and $G$ a finite group,
the equivalence relation reduces to that of conjugation in $G$.

% XXXXXXXXXXXXXXXXXXXXXXXXXXXX	EXAMPLE

\begin{example}
For the case of a point $Q = \{ p \}$ with the trivial action of a finite group $G$
(see Example \ref{expt}), the above defined equivalence relation reduces to
conjugation within the group, and $\tilde{Q}$ consists of one point for each
conjugacy class.  For an element $g \in G$, the point $(p, (g)_G)$ corresponding
to $(g)$ has the trivial action of $C(g)$, the centralizer of $g \in G$.
\end{example}

% XXXXXXXXXXXXXXXXXXXXXXXXXXXX	EXAMPLE

\begin{example}
If $Q$ is the $\Z_k$-teardrop (see Example \ref{extear}), then $\tilde{Q}$ has $k$
connected components.  The nontwisted sector is diffeomorphic to $Q$, while each
of the $k-1$ other components are points with the trivial action of $\Z_k$.
\end{example}

% XXXXXXXXXXXXXXXXXXXXXXXXXXXX	EXAMPLE

\begin{example}
For $Q = \R^3/\Z_2$ where $\Z_2$ acts via the antipodal map (see Example \ref{exrp2}),
$\tilde{Q}$ has two connected components, one diffeomorphic to $Q$ and one given
by a point with trivial $\Z_2$-action.
\end{example}

% XXXXXXXXXXXXXXXXXXXXXXXXXXXX	EXAMPLE

\begin{example}
Consider the case where $Q$ is the $\Z_k$-$\Z_l$-solid hollow football (see Example
\ref{exsolidhollow}) with $k$ and $l$ relatively prime.  Then $\tilde{Q}$
has $k + l - 1$ connected components.  The component corresponding to the identity
element (in both groups) is diffeomorphic to $Q$.  Each of the $k + l - 2$
other components is diffeomorphic to a closed interval $[0, 1]$ with the trivial action of
$\Z_k$ or $\Z_l$ (there are $k - 1$ components with trivial $\Z_k$-action and
$l - 1$ components with trivial $\Z_l$-action).

If $k = l$, then the equivalence relation on conjugacy classes defines
an isomorphism between the two $\Z_k$, and $Q$ is a global quotient.  Then there are exactly $k$
connected components of $\tilde{Q}$: one, corresponding to the identity in $\Z_k$, is
diffeomorphic to $Q$, while each of the others is given by $[-2, -1] \cup [1, 2]$
with trivial $\Z_k$-action (these two pieces correspond to the fixed-point set of
the nonidentity elements of $\Z_k$, which lie on the $z$-axis).

If $k$ and $l$ are not relatively prime, then letting $s:= \,\mbox{gcd}\,(k,l)$,
$Q$ can be expressed as a global quotient of the $\Z_{k/s}$-$\Z_{l/s}$-football by
$\Z_s$, and then the structure of $\tilde{Q}$ can be determined using both of
the above constructions.
\end{example}

% XXXXXXXXXXXXXXXXXXXXXXXXXXXX	EXAMPLE

\begin{example}
Suppose $Q$ is an orbifold with boundary homeomorphic to the closed 3-disk whose
singular set, located on the interior of the disk, is that given in Example
\ref{exfig8}.  Then there are six equivalence classes in $T$, each containing
exactly one element of $G_1 \cong D_6$.  The nontwisted sector is again diffeomorphic
to $Q$.  The twisted sector corresponding to the equivalence class of
$R_\pi^x$ is given by $S^1$ with a trivial
$\langle R_\pi^x \rangle \cong \Z_2$-action.  The twisted sectors corresponding
to $R_{\frac{2\pi}{3}}^z$ and $(R_{\frac{2\pi}{3}}^z)^2$ are both given by
$S^1$ with trivial $\langle R_{\frac{2\pi}{3}}^z \rangle \cong \Z_3$-action.
The twisted sectors corresponding to the elements $(R_\pi^x)(R_{\frac{2\pi}{3}}^z)$
and $(R_\pi^x)(R_{\frac{2\pi}{3}}^z)^2$ are both points with trivial $D_6$-action.
\end{example}

Now, suppose $Q$ admits an almost complex structure $J$.  As in the case without
boundary, $\tilde{Q}$ inherits an almost complex structure $\tilde{J}$ from $Q$
(note that this almost structure is on $\pi^\ast TQ$, which is not generally the
same as $T\tilde{Q}$.
For each $p \in Q$ and chart $\{ V_p, G_p, \pi_p \}$ at $p$, the almost complex
structure defines a homomorphism $\rho_p : G_p \rightarrow GL(n, \C)$ where $n$
is the dimension of $Q$ over $\C$.  For each $g \in G_p$, $\rho_p (g)$ can be
expressed as
\[
	\left[\begin{array}{cccc}
				e^{2\pi i m_{1, g}/m_g}	&	0	& \cdots	&	0
								\\\\
				0	&	e^{2\pi i m_{2, g}/m_g}	& 		&	\vdots
								\\\\
				0	&	0			& \ddots	&	\vdots
								\\\\
				0	&	0	&\cdots&	e^{2\pi i m_{n, g}/m_g}
	\end{array}\right] 
\]
with $m_g$ the order of $\rho_p(g)$ and $0 \leq m_{i, g} < m_g$.  Following \cite{chenruan}, define
\[
\begin{array}{rcccl}
	\iota	&:&	\tilde{Q}		&	\rightarrow	&	\Q	\\\\
		&:&	(p, (g)_{G_p})	&	\mapsto	&
					\sum\limits_{i = 1}^n \frac{m_{i, g}}{m_g} ,
\end{array}
\]
and then $\iota$ defines a map which is constant on each $\tilde{Q}_{(g)}$.
Let $\iota_{(g)} := \iota[(p, (g)_{G_p})]$ for any
$(p, (g)_{G_p}) \in \tilde{Q}_{(g)}$.  We refer to $\iota_{(g)}$ as the {\bf degree
shift number} of $\tilde{Q}_{(g)}$.

% XXXXXXXXXXXXXXXXXXXXXXXXXXXX  DEFINITION

\begin{definition}
Let $Q$ be an orbifold with boundary that admits an almost complex
structure $J$.  The {\bf relative orbifold cohomology groups} are
defined to be
\[
	H_{orb}^d (Q) = \bigoplus\limits_{(g) \in T} H^{d - 2\iota_{(g)}}
			(Q_{(g)}),
\]
where $H^{d - 2\iota_{(g)}} (Q_{(g)})$ is either the singular or
de Rham relative cohomology group of $Q_{(g)}$.

\end{definition}

In the sequel, we will examine analogs of the theorems of the previous chapter
with the characteristic classes taken to be elements of these cohomology groups.

% XXXXXXXXXXXXXXXXXXXXXXXXXXXXXXXXXXXXXXXXXXXXXXXXXXXXXXXXXXXXXXXXXXX
%			Section: GB Theorem
% XXXXXXXXXXXXXXXXXXXXXXXXXXXXXXXXXXXXXXXXXXXXXXXXXXXXXXXXXXXXXXXXXXX

\section{The Gauss-Bonnet and Poincar\'{e}-Hopf Theorems in Chen-Ruan Cohomology}

We begin with a lemma.

Let $Q$ be a closed orbifold of dimension $n$, and let $\rho : E \rightarrow Q$
be an orbifold vector bundle of rank $k$.  Suppose further that $E$ has the
following property:
for any chart $\{ V_p, G_p, \pi_p \}$ at $p$, and any subgroup $H$ of $G_p$,
the codimension of the fixed point set of $H$ in the fiber $E_p$ is equal
to the codimension of the fixed point set of $H$ in $V_p$ (or is zero
if the codimension of the fixed point set in $V_p$ is greater than the rank
$k$ of $E$).  This is the case, for instance, for the tangent and cotangent
bundles, their exterior powers, etc. (which can be verified using a metric
and the exponential map).  As $E$ is an orbifold, we may
apply the construction to form $\tilde{E}$.

% XXXXXXXXXXXXXXXXXXXXXXXXXXXX	LEMMA

\begin{lemma}

With $E$, $Q$ as above, $\tilde{E}$ is naturally
an orbifold vector bundle over $\tilde{Q}$, and the rank of $\tilde{E}$
on each connected component $\tilde{Q}_{(g)}$ of dimension $l$ is $k - n + l$.
In particular, the fiber of $\tilde{E}$ is zero on any $\tilde{Q}_{(g)}$
with codimension larger than $k$.

\end{lemma}

Note that the restriction on the group actions on $E$ is stronger than the requirement
that $E$ be a so-called good orbifold vector bundle (see \cite{ruangwt} for the definition;
our definition of orbifold vector bundle coincides with Ruan's definition of a good orbifold
vector bundle).  However, we will be applying this result to the tangent bundle only.  For general
(good) orbifold vector bundles, $\tilde{E}$ is still naturally a vector bundle over $\tilde{Q}$,
but the rank of $\tilde{E}$ over a connected component of $\tilde{Q}$ depends on the group action
on $E$.

\noindent       {\it Proof:}

First, we note that as the local groups and injections of $E$ are precisely those of $Q$, the set
$T$ is identical for both orbifolds.
Let $\tilde{\rho} : \tilde{E} \rightarrow \tilde{Q}$ be defined by $\tilde{\rho}(e, (g)) = (\rho(e), (g))$.
Then $\tilde{\rho}$ is certainly well-defined, as $e$ is fixed by a group element if $\rho(e)$ is.

Fix a point $(e, (g)) \in \tilde{E}$, where $e \in E$ and $(g) \in T$, and let
$p := \rho(e) \in Q$ denote the projection of $e$.  Then by the definition of $\tilde{E}$,
$(e, (g))$ is contained in a uniformizing set $\{ V_e^g, C(g), \pi_{e, g} \}$ induced by a
uniformizing set $\{ V_e, C(g), \pi_{e, g} \}$ of $e$.  Moreover, we can take this orbifold
chart for $E$ to be a uniformizing system $(V_p \times \R^k, G_p, \tilde{\pi}_p)$ of
the rank $k$ orbifold bundle induced by a system $\{ V_p, G_p, \pi_p \}$ near $p$ in $Q$.
Hence, $V_e = V_p \times \R^k$, $G_e = G_p$, and $\pi_e = \tilde{\pi}_p$.

By our definition of orbifold vector bundle, the kernel of the $G_p$ action on the $\R^k$ fiber
over a point $y \in V_p$ is the kernel of the $G_p$ action on $y$.  Hence, the fixed point set
$V_e^g$ of $g$ in $V_e = V_p \times \R^k$ is $V_p^g \times \R^{k - n + l}$ (or the zero
vector space if $k - n + l \leq 0$), where $l$ is
the dimension of $V_p^g$ (i.e. $k$ minus the codimension of $V_p^g$ in $V_p$).
If we define the map $\tilde{\rho}_2 : V_e^g = V_p^g \times \R^{k - n + l} \rightarrow V_p^g$
by projection onto the first factor, then for any $f = (y, v) \in V_e^g$ (say $\pi_{p}(y) = q \in Q$),
we have that
\[
\begin{array}{rcl}
	\pi_{p, g} \circ \tilde{\rho}_2 (f)
	&=&
	\pi_{p, g} (y)						\\\\
	&=&
	(q, (g))						\\\\
	&=&
	\tilde{\rho} ((q, (g)), v)				\\\\
	&=&
	\tilde{\rho} \circ \pi_{e, g} (f),
\end{array}
\]
and hence that $\pi_{p, g} \circ \tilde{\rho}_2 = \tilde{\rho} \circ \pi_{e, g}$.

Hence,
\[
	\{ V_e^g, C(g), \pi_{e, g} \} = \{ V_p^g \times \R^{k -n + l}, C(g), \pi_{e, g} \} 
\]
is a rank $k - n + l$ uniformizing system for a bundle over the uniformizing system
$\{ V_p^g, C(g), \pi_{p, g} \}$ about $(p, (g)) \in \tilde{Q}$.

We have that any point $(e, (g))$ of $\tilde{E}$ is contained in a bundle
uniformizing system over a uniformizing system of $(\rho(e), (g)) \in \tilde{Q}$
with projection $\tilde{\rho}$.  Given a compatible cover ${\mathcal U}$ of $\tilde{Q}$
(which can be taken to be induced by a compatible cover of $Q$) an injection
$\tilde{\lambda} : \{ V_p^g, C(g), \pi_{p, g} \} \rightarrow \{ V_q^h, C(h), \pi_{q, h} \}$
is always the restriction of an injection
$\lambda : \{ V_p, G_p, \pi_p \} \rightarrow \{ V_q, G_q, \pi_q \}$ (by \cite{satake2}, Lemma 1, there
is a bijection between elements of $G_q$ and such injections; hence, restricting to the subgroup
$C(h)$ decreases the number of injections).  Then the transition map
$\tilde{\phi} : V_p^g \rightarrow Aut(\R^{k-n+l})$ is simply the restriction of the corresponding
transition map $\phi : V_p \rightarrow Aut(\R^n)$ for $E$ to the fixed-point set $V_p^g$.
Therefore, these uniformizing systems patch together to give $\tilde{E}$ the structure
of an orbibundle over $\tilde{Q}$.

\bigskip \noindent Q.E.D. \bigskip

Applying the above argument to the tangent bundle shows that
$T(\tilde{Q}) = \widetilde{TQ}$; i.e. the tangent bundle of $\tilde{Q}$ is the collection
of twisted sectors of the tangent bundle of $Q$.  Similarly, the constructions of the
cotangent, exterior power, and tensor bundles commute with this construction.  Moreover,
any smooth section $\omega$ of the bundle $E$ naturally induces a smooth section
$\tilde{\omega}$ of the bundle $\tilde{E}$ via $\tilde{\omega} : (p, (g)) \mapsto (\omega(p), (g))$.

% XXXXXXXXXXXXXXXXXXXXXXXXXXXX	THEOREM

\begin{theorem}[The Second Gauss-Bonnet Theorem for Closed Orbifolds]
\label{gbcr}
Let $Q$ be a closed oriented orbifold of dimension $n$, and suppose $Q$ carries a connection $\omega$
with curvature $\Omega$.  Let $\tilde{\omega}$ denote the induced
connection $\omega$ on $\widetilde{TQ}$, and let $\tilde{\Omega}$
denote its curvature.  Let $E(\tilde{\Omega})$ denote the Euler curvature
form of $\tilde{\Omega}$, defined on $\tilde{Q}$, and then
\[
	\int\limits_{\tilde{Q}} E(\tilde{\Omega})	=	\chi (Q),
\]
the Euler characteristic of the underlying topological space of $Q$.  Moreover, if $Q$ is almost
complex, then $E(\tilde{\Omega})$ represents an element of the cohomology ring
$H_{orb}^\ast (Q)$ which is independent of the connection on $Q$.

\end{theorem}

\noindent	{\it Proof:}

We first note that on each $\tilde{Q}_{(g)}$,
$E(\tilde{\Omega})$ is a representative of the Euler class of $\tilde{Q}_{(g)}$,
so that, by the Gauss-Bonnet theorem for orbifolds,
\[
\begin{array}{rcl}
	\int\limits_{\tilde{Q}} E(\tilde{\Omega} )
	&=&
	\sum\limits_{(g) \in T} \;\;\;\; \int\limits_{\tilde{Q}_{(g)}} E(\tilde{\Omega})		\\\\
	&=&
	\sum\limits_{(g) \in T} \chi_{orb} (\tilde{Q}_{(g)}),
\end{array}
\]
where again the sum is over the set $T$ of equivalence classes of local group elements.
Note that $\tilde{Q}$ is closed, so that it has a finite number of connected components
(i.e. $T$ is finite).

Now, let $K$ be a simplicial decomposition of $Q$ such that for each simplex $\sigma \in K$,
the order of the isotropy group $G_p$ of $p$ is constant on the interior of $\sigma$ (see
\cite{moerdijk}).  Let $\tilde{K}$ be the simplicial decomposition of $\tilde{Q}$
induced by $K$, and for each $\sigma \in K$, denote
by $\sigma_{(g)}$ the corresponding simplex in $\tilde{K}$ which lies in
$\tilde{Q}_{(g)}$.  As $K$ is finite, let $\sigma^0, \sigma^1 , \ldots , \sigma^k$
be an enumeration of the simplices in $K$ so that
$\tilde{K} = \{ \sigma_{(g)}^i : 0 \leq i \leq k, (g) \in T \}$.

For each $\sigma_{(g)}^i$, let $p_{(g)}^i$ be a point on the interior of the simplex,
and let $h_{(g)}^i$ be an element of $G_{p_{(g)}^i}$ such that $h_{(g)}^i \in (g)$.
Note that the order $|(h_{(g)}^i)_{G_{p_{(g)}^i}}|$ of the conjugacy class of $h_{(g)}^i$ in
$G_{p_{(g)}^i}$, as well as the orders $|G_{p_{(g)}^i}|$ and $|C(h_{(g)}^i)|$,
are independent of the choices of $p_{(g)}^i$ and $h_{(g)}^i$.  Hence,
\[
\begin{array}{rcl}
	\sum\limits_{(g) \in T} \chi_{orb} (\tilde{Q}_{(g)})
	&=&
	\sum\limits_{(g) \in T} \;\;  \sum\limits_{i=0}^k (-1)^{\mbox{dim}\, \sigma_{(g)}^i} \frac{1}{| C(h_{(g)}^i) |} 
											\\\\
	&&\mbox{(if there is no $\sigma_{(g)}^i$ for a specific $(g)$ and $i$, then let the term be zero)}
											\\\\
	&=&
	\sum\limits_{(g) \in T} \;\; \sum\limits_{i=0}^k (-1)^{\mbox{dim}\, \sigma_{(g)}^i} \frac{|(h_{(g)}^i)|}{|G_{p_{(g)}^i}|} 
											\\\\
	&&\mbox{as $|G_{p_{(g)}^i}|=|(h_{(g)}^i)||C(h_{(g)}^i)|$}
											\\\\
	&=&
	\sum\limits_{i=0}^k\;\;\;\;\; \sum\limits_{(g) \in T : (g) \cap G_{p_{(g)}^i} \neq \emptyset}
		(-1)^{\mbox{dim}\, \sigma_{(g)}^i} \frac{|(h_{(g)}^i)|}{|G_{p_{(g)}^i}|}
											\\\\
	&=&
	\sum\limits_{i=0}^k (-1)^{\mbox{dim}\, \sigma_{(1)}^i} \frac{|G_{p_{(g)}^i}|}{|G_{p_{(g)}^i}|} 
											\\\\
	&=&
	\sum\limits_{i=0}^k (-1)^{\mbox{dim}\, \sigma_{(1)}^i} 
											\\\\
	&=&
	\chi (Q).
\end{array}
\]

To finish the proof, suppose $Q$ is almost complex, and note that for each $(g) \in T$, as
$\tilde{Q}_{(g)}$ is on its own an orbifold,
$E(\tilde{\Omega})$ is a representative of the Euler class of $\tilde{Q}_{(g)}$
in $H^\ast (\tilde{Q}_{(g)})$.  Denote this
class $e(g)$, and then $E(\tilde{\Omega})$ represents the element
\[
	\bigoplus\limits_{(g) \in T} e(g) \in H^\ast_{orb}(Q),
\]
which is an invariant of the connection.  Note that if $\tilde{Q}_{(g)}$ has dimension $d_{(g)}$, then
$e(g)$ is an element of $H_{orb}^{d_{(g)} + 2 \iota(g)} (Q)$.

\bigskip \noindent Q.E.D. \bigskip

Note that it is inessential that $\tilde{\Omega}$ be defined as being induced by a connection
on $Q$; the theorem holds if we begin with an arbitrary connection on $\tilde{Q}$.

In the case of an almost complex, reduced orbifold $Q$ of dimension $n$, the cohomology
group $H_{orb}^{n} (Q)$ is isomorphic to the de Rham group $H^{n}(Q)$.  Hence, the top part
of $E(\tilde{\Omega})$ is a representative of the Euler class of $Q$ with respect to
this isomorphism.  In the case that $Q$ is not reduced, if $i$ denotes the number of
elements of $T$ whose representatives act trivially, then
$H_{orb}^{n} (Q)$ is isomorphic to $H^{n} (Q) \oplus H^{n} (Q) \oplus \cdots \oplus H^{n} (Q)$
($i$ copies).  Then the top part of $E(\tilde{\Omega})$ is $i$ copies of the Euler curvature form.

% XXXXXXXXXXXXXXXXXXXXXXXXXXXX	DEFINITION

\begin{definition}

Let $Q$ be an almost complex, closed, oriented orbifold, $E \rightarrow Q$ a vector bundle, and
$\tilde{E} \rightarrow \tilde{Q}$ the induced bundle.  The {\bf orbifold Euler class}
$e_{orb}(E)$ is the cohomology class represented by $E(\tilde{\Omega})$ in $H_{orb}^{\ast}(Q)$
for some connection $\tilde{\omega}$ on $\tilde{Q}$ with curvature $\tilde{\Omega}$.

\end{definition}

% XXXXXXXXXXXXXXXXXXXXXXXXXXXX	COROLLARY

\begin{corollary}[The Second Poincar\'{e}-Hopf Theorem for Closed Orbifolds]
\label{phcr}
Let $X$ be a vector field on the closed orbifold $Q$ with a finite
number of zeros, and let $\tilde{X}$ be the induced vector field on $\tilde{Q}$.
Then
\[
	\mbox{ind} (\tilde{X}) = \chi (Q).
\]

\end{corollary}

\noindent	{\it Proof:}

This follows from the Poincar\'{e}-Hopf theorem for closed orbifolds \cite{satake2},
applied to each connected component of $\tilde{Q}$.  Note that, as vector fields must be tangent
to the singular set, a vector field with a finite number of zeros on $Q$ will induce a vector
field with a finite number of zeros on $\tilde{Q}$.

\bigskip \noindent Q.E.D. \bigskip

Again, it is inessential that we begin with a vector field on $Q$ and pull back to $\tilde{Q}$.

% XXXXXXXXXXXXXXXXXXXXXXXXXXXX	Subsection: Examples

\subsection{Examples}

% XXXXXXXXXXXXXXXXXXXXXXXXXXXX	EXAMPLE

\begin{example}

We start with the example of a single point $Q = \{ p \}$ with the trivial action of a finite group $G$.
In this case, the equivalence relation reduces to conjugation in the group.
Then $\tilde{Q} = \{ (p, (g)) : (g) \in T \}$,
and the degree shifting number $\iota_{(g)} = 0$ for each $(g) \in T$ (see \cite{chenruan}).

The contribution of each connected component $\{ (p, (g) ) \}$ of $\tilde{Q}$ to the orbifold
cohomology is in $H_{orb}^0 (Q)$, so that if $n$ is the number of conjugacy classes in $G$,
\[
	H_{orb}^d (Q) =
	\left\{	\begin{array}{ll}
		\R^n ,	&	d = 0,	\\\\
		0		&	d \neq 0
	\end{array} \right.
\]

The curvature form of each point is the function $\Omega(p, (g)) = \frac{1}{|(g)|} = \chi_{orb}(p, (g))$,
so that summing this value over each of the connected components of $\tilde{Q}$ gives
\[
\begin{array}{rcl}
	\sum\limits_{(g) \in T}	\frac{1}{|C(g)|}
	&=&
	\sum\limits_{(g) \in T} \frac{|(g)|}{|G|}				\\\\
	&=&
	1									\\\\
	&=&
	\chi (Q).
\end{array}
\]
\end{example}

% XXXXXXXXXXXXXXXXXXXXXXXXXXXX	EXAMPLE

\begin{example}
Let $Q$ denote the $\Z_k$-teardrop, and then
$\tilde{Q} = Q \sqcup \bigsqcup\limits_{i=1}^{k - 1} \{ (p, (i)) \}$.  The group $\Z_k$
acts trivially on each $(p,(i))$, so that the orbifold Euler characteristic of each
of these points is $\frac{1}{k}$.  Then $\chi_{orb} (Q) = \frac{k + 1}{k}$, so that
\[
\begin{array}{rcl}
	\int_{\tilde{Q}} E(\tilde{\Omega})
	&=&		\chi_{orb}(Q) + \sum\limits_{i=1}^{k-1} \chi_{orb}(p, (i))		\\\\
	&=&		\frac{k + 1}{k} + (k - 1)\frac{1}{k}					\\\\
	&=&		2											\\\\
	&=&		\chi (Q) .
\end{array}
\]
\end{example}

% XXXXXXXXXXXXXXXXXXXXXXXXXXXXXXXXXXXXXXXXXXXXXXXXXXXXXXXXXXXXXXXXXXX
%			Section: Reults for boundary case
% XXXXXXXXXXXXXXXXXXXXXXXXXXXXXXXXXXXXXXXXXXXXXXXXXXXXXXXXXXXXXXXXXXX

\section{The Case With Boundary}

We return to the case of an orbifold with boundary.  A modification of the proof
of \ref{gbcr} shows:

% XXXXXXXXXXXXXXXXXXXXXXXXXXXX	PROPOSITION

\begin{theorem}[The Second Gauss-Bonnet Theorem for Orbifolds with Boundary]
\label{gbbdcr}
Let $Q$ be a closed orbifold of dimension $n$ with boundary $M$,
and suppose $Q$ carries a connection $\omega$
with curvature $\Omega$.  Let, $\tilde{\omega}$, $\tilde{\Omega}$, etc. be defined
as in the previous section, and then
\[
	\int\limits_{\tilde{Q}} E(\tilde{\Omega})	=	\chi^\prime (Q)
							- \frac{1}{2} \chi(M).
\]
If $Q$ is almost complex, then $E(\tilde{\Omega})$ represents an element of the
cohomology ring $H_{orb}^\ast (Q)$ which is independent of the connection on $Q$.

\end{theorem}

Note that $\chi^\prime (Q) = \chi(Q, M)$.

\noindent	{\it Proof:}

Again, by the Gauss-Bonnet theorem for orbifolds with boundary,
\[
\begin{array}{rcl}
	\int\limits_{\tilde{Q}} E(\tilde{\Omega} )
	&=&
	\sum\limits_{(g) \in T} \;\;\;\; \int\limits_{\tilde{Q}_{(g)}} E(\tilde{\Omega})		\\\\
	&=&
	\sum\limits_{(g) \in T} \chi^\prime_{orb} (\tilde{Q}_{(g)})
		- \frac{1}{2} \chi_{orb} (M_{(g)}) .
\end{array}
\]
Using a simplicial decomposition of $Q$ as above and the same counting argument,
we have
\[
	\sum\limits_{(g) \in T} \chi^\prime_{orb} (\tilde{Q}_{(g)})
		- \frac{1}{2} \chi_{orb} (M_{(g)}) 
	=
	\chi^\prime  (\tilde{Q}) - \frac{1}{2} \chi  (M) .
\]

\bigskip \noindent Q.E.D. \bigskip

Clearly, in the case where the dimension $n$ of $Q$ is even, this formula becomes
\[
	\int\limits_{\tilde{Q}} E(\tilde{\Omega})	=	\chi^\prime (Q),
\]
and in the case where $n$ is odd,
\[
	\frac{1}{2} \chi (M)	=	\chi^\prime (Q).
\]

We now return to the proof of Proposition \ref{phdependent}.  Let $\tilde{\Psi}$ be defined
in the natural way by taking the sum of $\Psi$ on each connected component
of $\tilde{Q}$.  Then the relation $- d\tilde{\Psi} = \rho^\ast E(\tilde{\Omega})$ is
immediate, as it is true on each connected component.  Again, let $\tilde{X}$
denote the extension of the vector field $X$ to $\tilde{Q}$.  We modify the proof
of Proposition \ref{phdependent} by applying in the first step Theorem \ref{gbbdcr} instead
of the Theorem \ref{gbb1} as follows.

% XXXXXXXXXXXXXXXXXXXXXXXXXXXX	even case

If the dimension $n = 2m$ of $Q$ is even, then
\[
\begin{array}{rcl}
	\chi^\prime (Q)
	&=&
	\int\limits_Q E(\tilde{\Omega})						\\\\
	&&\mbox{(by Theorem \ref{gbbdcr})}					\\\\
	&=&
	\lim\limits_{r \to 0^+} \int\limits_{\tilde{Q} \backslash B_r(p)}
		\alpha^\ast \rho^\ast (E(\tilde{\Omega}))		\\\\
	&=&
	-\lim\limits_{r \to 0^+} \int\limits_{\tilde{Q} \backslash B_r(p)}
		d\alpha^\ast (\tilde{\Psi})					\\\\
	&=&
	\lim\limits_{r \to 0^+}
	\int\limits_{\partial B_r(p)} \alpha^\ast (\tilde{\Psi})
	- \int\limits_{\tilde{M}} \alpha^\ast(\tilde{\Psi})
										\\\\
	&=&
	\mbox{ind}(\tilde{X}) - \int\limits_{\tilde{M}} \alpha^\ast(\tilde{\Psi}),
\end{array}
\]
and hence
\[
	\mbox{ind}(\tilde{X})
	=
	\chi^\prime (Q) + \int\limits_{\tilde{M}} \alpha^\ast(\tilde{\Psi}).
\]
Note that the singular points $p$ are taken to be those of
$\tilde{X}$, and hence the $B_r(p)$ contains balls each singular
point of $\tilde{X}$.

% XXXXXXXXXXXXXXXXXXXXXXXXXXXX	odd case

Making the identical modification to the proof in the case that
$n = 2m + 1$ is odd, we obtain
\[
	\chi^\prime (Q) - \frac{1}{2} \chi (M)
	=
	\mbox{ind}(\tilde{X}) - \int_{\tilde{M}} \alpha^\ast(\tilde{\Psi}),
\]
and hence
\[
	\mbox{ind}(\tilde{X})
	=
	\chi^\prime (Q) - \frac{1}{2} \chi (M)
		+ \int_M \alpha^\ast(\tilde{\Psi}) .
\]

Now, note that in the case that $Q$ admits a complex structure,
as the cohomology class of $\Psi$ is independent of
the connection chosen, the cohomology class $\tilde{\Upsilon}$
of $\tilde{\Psi}$ in $H_{orb}^\ast(STQ_{|M})$ is similarly
independent.  In fact, it is clear that we can define $\tilde{\Upsilon}$
to be the sum of the cohomology classes of the forms $\Psi$ defined
on each connected component of $\tilde{Q}$ from the connection, and then
$\tilde{\Psi}$ would be a representative of the cohomology class
$\tilde{\Upsilon}$.  Hence, we have proven

% XXXXXXXXXXXXXXXXXXXXXXXXXXXX	THEOREM

\begin{theorem}
[The Second Poincar\'{e}-Hopf Theorem for Orbifolds with Boundary]

\label{pccr}
Let $Q$ be a compact oriented orbifold with  boundary $M$, and suppose
$Q$ admits an almost complex structure.  Let
$X$ be vector field on $Q$ which has a finite number of singularities,
all of which occuring on the interior of $Q$.  Then with $\tilde{Q}$,
$\tilde{X}$, etc. defined as above, we have
\[
	\mbox{ind}(\tilde{X})		=
		\chi^\prime (Q) +	\tilde{\alpha}^\ast \tilde{\Upsilon} ([\tilde{M}]) .
\]
\end{theorem}

Here, $\tilde{\alpha}^\ast \tilde{\Upsilon} ([\tilde{M}])$ refers to the integral of any form
representing the cohomology class $\tilde{\Upsilon}$ over the orbifold
$\tilde{M}$; we have chosen to use this notation to emphasize the fact
that the value of this integral is independent of the particular
representative of $\tilde{\Upsilon}$ chosen.

Note that in the case that $Q$ is a smooth manifold without boundary, both
Theorem \ref{pcb1} and Corollary \ref{pccr} reduce to the classical
Poincar\'{e}-Hopf Theorem.  Hence, both can be considered to be
generalizations of this theorem to orbifolds with boundary, in the
spirit of \cite{sha}.

%% file: thesisappA.tex
% Thesis Appendix 4, Christopher Seaton, seatonc@colorado.edu
% XXXXXXXXXXXXXXXXXXXXXXXXXXXXXXXXXXXXXXXXXXXXXXXXXXXXXXXXXXXXXXXXXXX
% XXXXXXXXXXXXXXXXXXXXXXXXXXXXXXXXXXXXXXXXXXXXXXXXXXXXXXXXXXXXXXXXXXX
%                       Appendix A: Orb Euler Char.
% XXXXXXXXXXXXXXXXXXXXXXXXXXXXXXXXXXXXXXXXXXXXXXXXXXXXXXXXXXXXXXXXXXX
% XXXXXXXXXXXXXXXXXXXXXXXXXXXXXXXXXXXXXXXXXXXXXXXXXXXXXXXXXXXXXXXXXXX
%

\chapter{A Note on Orbifold Euler Characteristics}
\label{appa}

This note is intended to clarify the various definitions of Euler Characteristics
for that have been given for orbifolds and preclude any confusion which may arise.
Throughout, let $Q$ be a compact oriented orbifold and let $X_Q$ denote
its underlying topological space.

The first orbifold Euler characteristic was defined by Satake in \cite{satake2},
and was originally denoted $\chi_V (Q)$.  This is the orbifold Euler
characteristic used in this paper, and which we have denoted $\chi_{orb}(Q)$.
It was defined as the index, in the
orbifold sense, of a vector field on $Q$ with isoloated singularities, or
equivalently as the integral, in the orbifold sense, of the Euler curvature
form $E(\Omega)$, defined in terms of an orbifold connection on $Q$.
In general, it is a rational number.
If $e(Q)$ denotes the cohomology class of $E(\Omega)$ and $e(X_Q)$
the usual Euler class of the underlying space of $Q$, then using
the isomorphism between the de Rahm cohomology of $Q$ as an orbifold
and the singular cohomology of $X_Q$, it is easy to see that
\[
	\frac{ \chi_{orb} (Q) }{\chi(X_Q)} e(Q) = e(X_Q) .
\]
In many ways, this Euler characteristic plays the role of the usual
Euler characteristic of a manifold; as was discussed in Chapter \ref{chgbph},
it can be defined in terms of suitable triangulation of the orbifold.

In the case that $Q = M/G$ is the quotient of a manifold by a finite
group, another orbifold Euler number was defined in \cite{dixon} as
\[
	\chi (M, G) := \frac{1}{|G|} \sum\limits_{g, h \in G :gh = hg}
		\chi (M^g \cap M^h),
\]
where $M^g$ denotes the fixed-point set of $g$ (see also \cite{atiyahsegal},
where this is shown to be the Euler characteristic of the equivariant
$K$-theory $K_G(M)$).  In \cite{roan}, this definition is generalized to the
case of a general complex orbifold with groups acting as subgroups of
$SL_n(\C)$ as follows:  If $e(p)$ denotes for each $p \in Q$
the number of conjugacy classes in the isotropy subgroup $I_p$,
then define
\[
	\chi^o (Q) := \sum\limits_{k = 1}^\infty k \chi(e^{-1}(k)) .
\]
Note that the sum on the right has finitely many non-zero terms.
Roan shows that this definition coincides with that of Dixon {\it et. al.}
\cite{dixon}
in the case of a global quotient, and moreover that it does not
coincide with the usual Euler characteristic of the underlying space of
$Q$ unless $Q$ is nonsingular.

Similarly, in \cite{ademruan}, the formula for $\chi(M, G)$ is genearlized to general
orbifolds by defining
\[
	\chi_{orb}(Q) := \:\mbox{dim}\: K_{orb}^0 (Q) \otimes \Q -
		\:\mbox{dim}\: K_{orb}^1 (Q) \otimes \Q .
\]
Here, $K_{orb}(Q)$ denotes the {\bf orbifold $K$-theory} of $Q$, which is defined
using complex orbifold vector bundles over $Q$ and generalizes the equivariant
$K$-theory of $M/G$.

We note that in the case of a reduced orbifold,
these two invariants $\chi^o(Q)$ and $\chi_{orb}(Q)$ coincide.
Indeed, using the notation of Chapter \ref{chcr},
for each $p \in Q$, $e(p)$ coincides with the number of distinct
$(g) \in T$ such that $(p, (g)) \in \tilde{Q}$.  Hence, for each positive integer
$k$, the preimage $\pi^{-1} e^{-1} (k)$ of $e^{-1}(k)$ in $\tilde{Q}$ via the
projection $\pi : \tilde{Q} \rightarrow Q$ is a $k$-fold disjoint covering.
Therefore, we have that
\[
\begin{array}{rcl}
	\chi^o (Q)	&=&	\sum\limits_{k = 1}^\infty k \chi(e^{-1}(k))		\\\\
			&=&	\sum\limits_{(g) \in T} \chi (\tilde{Q}_{(g)})		\\\\
			&=&	\chi_{orb} (Q) .
\end{array}
\]
The last equality follows from the decomposition theorem for orbifold $K$-theory of Adem-Ruan
(\cite{ademruan} Theorem 5.1, Corollary 5.6).

%% file: thesisappB.tex
% Thesis Appendix 4, Christopher Seaton, seatonc@colorado.edu
% XXXXXXXXXXXXXXXXXXXXXXXXXXXXXXXXXXXXXXXXXXXXXXXXXXXXXXXXXXXXXXXXXXX
% XXXXXXXXXXXXXXXXXXXXXXXXXXXXXXXXXXXXXXXXXXXXXXXXXXXXXXXXXXXXXXXXXXX
%                       Appendix B: Index Theorems
% XXXXXXXXXXXXXXXXXXXXXXXXXXXXXXXXXXXXXXXXXXXXXXXXXXXXXXXXXXXXXXXXXXX
% XXXXXXXXXXXXXXXXXXXXXXXXXXXXXXXXXXXXXXXXXXXXXXXXXXXXXXXXXXXXXXXXXXX
%

\chapter{Relation to the Index Theorem for Orbifolds}
\label{appb}

Let $Q$ be a compact orbifold that admits an almost complex structure,
$TQ$ its (orbifold) tangent bundle. Let $E$ and $F$ be two orbifold vector
bundles over $Q$, and $P$ an elliptic pseudodifferential operator
$C^\infty (E) \rightarrow C^\infty (F)$ from the smooth sections of
$E$ to those of $F$.  Then if $\rho : T^\ast Q \rightarrow Q$ denotes
the projection of the tangent bundle, the symbol $\sigma$ of $P$ gives
an isomorphism $\sigma : \rho^\ast E \rightarrow \rho^\ast F$ off of
the zero-section of $T^\ast Q$, and hence an element of $K_{orb}(T^\ast Q)$.

Let $u \in K_{orb}(T^\ast Q)$ denote such a class in the orbifold $K$-theory
of $T^\ast Q \cong TQ$.  Let $\tilde{Q}, \tilde{Q}_{(g)}$, etc. be defined as in
Chapter \ref{chcr}, and let $\tilde{E}, \tilde{F}$ denote the induced
bundles on $\tilde{Q}$.  Let $\tilde{E}_{(g)}$ denote the induced
bundle restricted to the connected component $\tilde{Q}_{(g)}$
of $\tilde{Q}$, and let $u_{(g)}$ denote the class of
$\tilde{E}_{(g)} - \tilde{F}_{(g)}$ in $K_{orb} (\widetilde{TQ}_{(g)})
= K_{orb} (T \tilde{Q}_{(g)} )$ with the symbol given by the restriction of $P$.
With this setup, Kawasaki \cite{kawasaki2} showed that the orbifold index
\[
	\mbox{ind}_{orb} (P)	:=	\mbox{dim} [\mbox{ker}(P)]
				-	\mbox{dim} [\mbox{coker}(P)]
\]
is given by
\[
	\mbox{ind}_{orb} (u) =
	(-1)^{\mbox{dim}\: Q} \langle \mbox{ch}(u) \cup {\mathcal T}(Q), [TQ] \rangle
		+ \sum\limits_{(1) \neq (g) \in T} \frac{(-1)^{\mbox{dim}\: \tilde{Q}_{(g)}}}{|C(g)|}
		\langle \mbox{ch}(u_{(g)}) \cup {\mathcal T}(\tilde{Q}_{(g)}), [T\tilde{Q}_{(g)}] \rangle
\]
(see also \cite{farsi}).
Here, ch$(u)$ denotes the Chern character of $u$ (ch$(E) - $ch$(F)$), and
${\mathcal T}$ the Todd class of te (complexified) tangent bundle.

We we may define the {\bf orbifold Chern character} $\mbox{ch}_{orb}(E)$
of a complex vector bundle $E$ and the {\bf orbifold Todd class} ${\mathcal T}_{orb}(Q)$
in an analogous manner to the orbifold Euler class: to a bundle
$E$ with connection $\omega$ we associate the induced bundle $\tilde{E}$ over
$\tilde{Q}$ with connection $\tilde{\omega}$, and then construct the the Chern
character and Todd class from the curvature of this connection in the usual
manner (the bundle is taken to be the tangent bundle in the case of the Todd
class).  It is clear that these forms represent elements of $H^\ast(Q)$
that are independent of the connection, and that they can be given as a sum
of ordinary characteristic classes of the $\tilde{Q}_{(g)}$.  In other words,
we have
\[
	\mbox{ch}_{orb}(\tilde{E} - \tilde{F})
		= \sum\limits_{(g) \in T}
		\: \mbox{ch}(\tilde{E}_{(g)} - \tilde{F}_{(g)})
\]
and
\[
	{\mathcal T}_{orb}(Q)
		= \sum\limits_{(g) \in T}
		\: {\mathcal T}(\tilde{Q}_{(g)}).
\]
We note that the factors $\frac{1}{|C(g)|}$ appear in the above formula because
Kawasaki worked with a definition of orbifold which required that it was reduced.
Hence, his $\tilde{Q}_{(g)}$ correspond to our $(\tilde{Q}_{(g)})_{red}$, and
the corresponding factors are built into our definition of the integral.

With this, by reversing the orientation of the odd-dimensional components of
$\tilde{Q}$, Kawasaki's index formula can be written as 
\[
	\mbox{ind}_{orb} (u) =
	\langle \mbox{ch}_{orb}(u) \cup {\mathcal T}_{orb}(Q), [\widetilde{T Q}] \rangle .
\]

We note that the cup product here is {\bf not} the orbifold cup product defined in
\cite{chenruan}, but rather the pointwise product induced by the wedge product.
In particular, forms on $\tilde{Q}_{(g)}$ cupped with forms on $\tilde{Q}_{(h)}$
yield $0$ when $(g) \not\sim (h)$.

%% file: thesismaster.bbl
\begin{thebibliography}{10}

\bibitem{ademruan}
A.~Adem and Y.~Ruan.
\newblock Twisted orbifold ${K}$-theory.
\newblock math.AT/0107168, 2001.

\bibitem{borzellino}
J.~Borzellino.
\newblock {\em {R}iemannian Geometry of Orbifolds}.
\newblock PhD thesis, UCLA, 1992.

\bibitem{chenruan}
W.~Chen and Y.~Ruan.
\newblock A new cohomology theory for orbifold.
\newblock math.AG/0004129, 2001.

\bibitem{chern1}
S.S. Chern.
\newblock A simple intrinsic proof of the {G}auss-{B}onnet formula for closed
  {R}iemannian manifolds.
\newblock {\em Ann. of Math.}, 45:747--752, 1944.

\bibitem{chern2}
S.S. Chern.
\newblock On the curvatura integra in a {R}iemannian manifold.
\newblock {\em Ann. of Math.}, 99:48--69, 1974.

\bibitem{dixon}
L.~Dixon, J.~Harvey, C.~Vafa, and E.~Witten.
\newblock Strings on orbifolds.
\newblock {\em Nucl. Phys. B}, 261:678--686, 1985.

\bibitem{farsi}
C.~Farsi.
\newblock ${K}$-theoretical index theorems for good orbifolds.
\newblock {\em Proceedings of the AMS}, 3:769--773, 1992.

\bibitem{gp}
V.~Guillemin and A.~Pollack.
\newblock {\em Differential Topology}.
\newblock Prentice-Hall, Inc., Englewood Cliffs, New Jersey, 1965.

\bibitem{kawasaki1}
T.~Kawasaki.
\newblock The signature theorem for ${V}$-manifolds.
\newblock {\em Topology}, 17:75--83, 1977.

\bibitem{kawasaki2}
T.~Kawasaki.
\newblock The index of elliptic operators over ${V}$-manifolds.
\newblock {\em Nagoya Math. J.}, 84:135--157, 1981.

\bibitem{liang}
C.~Liang.
\newblock Vector fields on ${V}$-manifolds, and locally free ${G} \times
  {{\R}}^l$-actions on manifolds.
\newblock {\em Indiana University Mathematics Journal}, 27:349--352, 1978.

\bibitem{atiyahsegal}
G.~Segal M.F.~Atiyah.
\newblock On equivariant {E}uler characteristics.
\newblock {\em J. Geom. Phys.}, 4:671--677, 1989.

\bibitem{milnor2}
J.~Milnor.
\newblock {\em Topology from the Differentiable Viewpoint}.
\newblock The University Press of Virginia, Charlottesville, 1965.

\bibitem{milnor}
J.~Milnor and J.~Stasheff.
\newblock {\em Annals of Mathematical Studies 76: Characteristic Classes}.
\newblock Princeton University Press and University of Tokyo Press, Princeton,
  New Jersey, 1974.

\bibitem{moerdijk}
I.~Moerdijk and D.~A. Pronk.
\newblock Simplicial cohomology of orbifolds.
\newblock math.alg/9708021, 1997.

\bibitem{roan}
S.~S. Roan.
\newblock Minimal resolutions of {G}orenstein orbifolds in dimension three.
\newblock {\em Topology}, 35:489--509, 1996.

\bibitem{ruangwt}
Y.~Ruan.
\newblock Orbifold {G}romov-{W}itten theory.
\newblock math.AG/0103156, 2001.

\bibitem{ruansgt}
Y.~Ruan.
\newblock Stringy geometry and topology of orbifolds.
\newblock math.AG/0011149, 2001.

\bibitem{satake1}
I.~Satake.
\newblock On a generalization of the notion of manifold.
\newblock {\em Proc. Nat. Acad. Sci. USA}, 42:359--363, 1956.

\bibitem{satake2}
I.~Satake.
\newblock The {G}auss-{B}onnet theorem for ${V}$-manifolds.
\newblock {\em Journ. Math. Soc. Japan}, 9:464--492, 1957.

\bibitem{schwerdtfeger}
B.~E. Schwerdtfeger.
\newblock {\em On Explicit Furmulas for the Differential Forms in the Thom
  Class}.
\newblock PhD thesis, University of Bonn, 1979.

\bibitem{sha}
J.~P. Sha.
\newblock A secondary {C}hern-{E}uler class.
\newblock {\em Ann. of Math.}, 159:1151--1158, 1999.

\bibitem{stanhope}
E.~A. Stanhope.
\newblock {\em Hearing Orbifold Topology}.
\newblock PhD thesis, Dartmouth College, 2002.

\bibitem{thurston}
W.~Thurston.
\newblock {\em The Geometry and Topology of 3-Manifolds}.
\newblock Lecture Notes, Princeton University Math Dept., Princeton, New
  Jersey, 1978.

\end{thebibliography}
